\documentclass[12pt]{amsart}
\makeatletter
\let\atopwithdelims\@@atopwithdelims
\makeatother
\def\Da{{\Delta}}

\def\cC{{\mathcal C}}

\def\ZZ{{\mathbb Z}}
\def\CC{{\mathbb C}}
\def\RR{{\mathbb R}}
\def\NN{{\mathbb N}}
\def\Lat{{\mathcal Lat}}

\def\al{\alpha}
\def\lie{\ell ie}

\def\int#1{\mathop{#1}\nolimits_{}^{\circ}}

\def\sta,{\mbox{\rm star}_\Delta}

\def\id{{\mbox{\rm id}}}
\def\#{\sharp}

\def\Cut{{\rm Cut}}

\def\log{{\rm log}}
\def\exp{{\rm exp}}
\newcommand{\lb}{\left[}
\newcommand{\rb}{\right]}
\newcommand{\lp}{\left\{}
\newcommand{\rp}{\right\}}

\newcommand{\de}{\Delta}
\newcommand{\dee}[2]{\de_{#1}^{#2}}
\newcommand{\ra}{\rightarrow}
\newcommand{\ov}[1]{\overline{#1}}
\newcommand{\ti}[1]{\tilde{#1}}
\newcommand{\mb}[1]{\mu_{#1}(\hat{0},\hat{1})}
\newcommand{\ho}{\hat{0}}
\newcommand{\me}[2]{\mu_{#1}(\ho,#2)}
\newcommand{\meg}{\me{M^g}{F}}
\newcommand{\mega}{\me{L^{g^\ast}}{G}}
\newcommand{\snk}{\Sigma_{n,k}}
\newcommand{\osnk}{\ov{\Sigma}_{n,k}}
\newtheorem{theorem}{Theorem}[section]
\newtheorem{lemma}[theorem]{Lemma}
\newtheorem{proposition}[theorem]{Proposition}
\newtheorem{corollary}[theorem]{Corollary}

\newenvironment{Proof}{\smallskip\noindent\mbox{\sf Proof: } }{\qed\medskip}

\newenvironment{definition}{\smallskip\noindent\mbox{\sf Definition: } }{\medskip}

\begin{document}
\title{Complexes of not $i$-Connected Graphs}

\author[E. Babson]{Eric Babson}
\thanks{Babson, Bj\"orner, Linusson and Welker
were partially supported by MSRI.
Babson was
supported by a National Science Foundation postdoctoral
fellowship.
Linusson was supported by a Swedish Natural Sciences Research Council 
(NFR) postdoctoral fellowship.
Welker was
supported by Deutsche Forschungsgemeinschaft (DFG).
Research at MSRI is supported in part by NSF grant DMS-9022140.
} 
\address{\hskip-\parindent Eric Babson,
Mathematical Sciences Research
Institute, 1000 Centennial Drive, Berkeley, CA 94720-5070, USA.}
\email{babson@msri.org}
\author[A. Bj\"orner]{Anders Bj\"orner}
\address{\hskip-\parindent Anders Bj\"orner, Department of Mathematics, Royal Institute of Technology,
S-100 44 Stockholm, Sweden}
\email{bjorner@math.kth.se}
\author[S. Linusson]{Svante Linusson}
\address{\hskip-\parindent Svante Linusson, Department of Mathematics, Stockholms Universitet,
S-106 91 Stockholm, Sweden}
\email{linusson@matematik.su.se}
\author[J. Shareshian]{John Shareshian}
\address{\hskip-\parindent John Shareshian, Mathematical Sciences Research 
Institute, 1000 Centennial Drive, Berkeley, CA 94720-5070, USA.}
\email{shareshi@msri.org}
\author[V. Welker]{Volkmar Welker}
\address{\hskip-\parindent Volkmar Welker, Fachbereich 6, Mathematik, Universit\"at GH-Essen, D-45117 Essen, Germany}
\email{welker@exp-math.uni-essen.de}



\date{}

\begin{abstract}
Complexes of (not) connected graphs, hypergraphs and their
homology appear in the construction of knot invariants given
by V. Vassiliev \cite{Vas93-1,Vas94,Vas97}. In this paper we 
study the complexes of not $i$-connected $k$-hypergraphs on $n$ vertices. 
We show that the complex of not $2$-connected  
graphs has the homotopy type of a wedge of $(n-2)!$ spheres of 
dimension $2n-5$. This answers one of the questions raised by 
Vassiliev \cite{Vas97} in connection with knot invariants. For this 
case the $S_n$-action on the homology of the complex is also determined. 
For complexes of not $2$-connected $k$-hypergraphs we provide a formula for the 
generating function of the Euler characteristic, and we introduce certain
lattices of graphs that encode their topology. We also present 
partial results for some other cases. In particular, we show that 
the complex of not $(n-2)$-connected graphs is Alexander dual to the 
complex of partial matchings of the complete graph. 
For not $(n-3)$-connected graphs we provide a formula for the 
generating function of the Euler characteristic. \end{abstract}

\maketitle

\section{Introduction}

In this paper we study the homotopy type and homology of simplicial complexes 
whose simplices are
the edge sets of not $i$-connected graphs and hypergraphs on $n$ vertices. 
The case $i=1$ is already well understood (see Proposition \ref{discon}),
and here we begin the examination of the topological structure of
such complexes for $i \geq 2$.

Although our point of view is mainly combinatorial,
our original motivation for studying these complexes 
comes from the theory of Vassiliev invariants in knot theory.
By determining the homotopy type of the complex of not 
$2$-connected graphs on $n$ vertices we answer a question posed by 
V. Vassiliev in \cite{Vas97}, where he presents a new approach 
to Vassiliev knot invariants using a filtration of the simplicial 
resolution of the space of not-knots as in \cite{Vas94}. More precisely,
he studies the space $\Sigma$ of maps $f : S^1 \rightarrow \RR^3$ such 
that $f(S^1)$ has multiple points or
cusps. The simplicial resolution $\widetilde{\Sigma}$ of
$\Sigma$ is obtained roughly speaking as follows:
singular knots are resolved by blowing up
each $r$-fold self-intersection to an $({r\choose 2}-1)$-simplex,
and similarly for the set of cusps.
A suitable filtration (see \cite{Vas97}) 
of $\widetilde{\Sigma}$, combinatorially defined in terms of these
simplices, gives rise to a spectral sequence that 
contains the homology of the complex of not $2$-connected graphs
on $n$ vertices as a basic ingredient. 

Our work continues the
already fruitful interaction between the theory of Vassiliev invariants 
and questions in topological and homological combinatorics of 
graph complexes (see \cite{Vas93-1}). The study of
complexes of not $i$-connected graphs has intriguing
combinatorial and algebraic aspects as well. For example, such aspects become 
apparent when considering the complex of not $(n-2)$-connected graphs
on $n$ vertices. In Section \ref{match} this complex
is shown to be Alexander dual to the complex of partial matchings 
of the complete graph on $n$ vertices. These matching complexes, along
with complexes of partial matchings of bipartite graphs, have previously
been studied for other reasons, see \cite{BLVZ92}. In each case for which 
we calculate
the Betti numbers, we detect nontrivial homology. For $(n-3)$-connected
graphs (see Section \ref{notn3}) and for most complexes of not 
$2$-connected hypergraphs (see Section \ref{notconnhyper})
we have been unable to compute the Betti numbers explicitly, but we do
determine the generating function of their reduced Euler characteristics.  The 
homology is seen to be nontrivial in almost all of these cases. 

Surprisingly, these non-vanishing phenomena
are suggested by a result motivated by a conjecture in complexity theory.
The conjecture states that complexes of graphs on $n$ vertices
having some non-trivial monotone graph property -- like being not 
$i$-connected -- 
are evasive (see for example \cite{KahSakStu84}). Kahn, Saks 
\& Sturtevant \cite{KahSakStu84} showed that non-evasive complexes are 
contractible. In many naturally arising cases, including those examined
here, the converse is true 
and evasive complexes in fact have non-vanishing 
reduced Euler characteristics. 

In Section \ref{rep} we study the action of the symmetric group on the
complex of not $2$-connected graphs induced by its natural action on the
vertices. This action induces a representation of $S_n$ on the homology
groups of the complex, which we determine. Using the representation, we
deduce upper bounds on the number of Vassiliev invariants of a given
bi-order.  This representation
coincides with a recently well studied representation which appears in 
the work of Robinson \& Whitehouse \cite{RobWhi94,Whi94-2},
Kontsevich \cite{Kon93}, Getzler \& Kapranov
\cite{GetKap95}, Mathieu \cite{Mat96}, Hanlon \& Stanley 
\cite{Han96,HanSta95} and Sundaram \cite{Sun96}. 

\medskip
\noindent {\sf Acknowledgment:} We are grateful to V. Vassiliev for inspiring 
discussions and hints, which sparked our interest and initiated this research.
All computer calculations presented in this paper were performed 
using a Mathematica$^{\text{\copyright}}$ package designed by Vic Reiner
and a C-Program by Frank Heckenbach.

\section{Preliminaries} \label{prel}

We now introduce the basic concepts used in this paper.
By a graph $G = (V(G),E(G))$ we mean a loopless graph without 
multiple edges on the vertex set $V(G)$ and with edge set $E(G) 
\subseteq \binom{V(G)}{2}$.  Our standard vertex set
will be the set $[n]:=\{1,2,\dots,n\}$.
A graph $G$ is called {\it connected} 
if for any two distinct vertices $v, v' \in V(G)$ there is a 
path from $v$ to $v'$ in $G$, that is, a sequence of edges 
$\{ v_1, v_2 \}$, $\{ v_2, v_3 \}$, $\ldots$, 
$\{v_{l-1},v_l\}$ $\in E(G)$ such that $v = v_1$ and $v' = v_{l}$. 
Such a path will sometimes be denoted by $v_1,v_2,\ldots,v_l$.  The
{\it size} of a graph $G$ is $|V(G)|$.

A graph $G$ is called {\it $i$-connected}, for a number $i$
such that $0<i<|V(G)|$,
if for any $j$ vertices $v_1, \ldots, v_j \in V(G)$, $j<i,$ the graph
$G'$ that is obtained from $G$ by deleting the vertices $v_1, 
\ldots, v_j$ and their adjacent edges is connected.  
Equivalently, $G$ is $i$-connected if and only if for every pair
$v, v'$ of not adjacent vertices there are at least $i$ paths
from $v$ to $v'$ that are pairwise vertex disjoint except at their endpoints.

A graph with at least $i+1$ vertices which
is not $i$-connected is also called {\it (i-1)-separable}, and a
$1$-separable (that is, not $2$-connected) graph will often be 
called just separable.
Of course, if $G = (V(G),E(G))$ is a graph that is not $i$-connected 
for some $i \geq 1$ then for any subset $E' \subseteq E(G)$ the 
graph $G' = (V(G),E')$ on the same vertex set is not $i$-connected 
either. Hence if we fix an $n$-element vertex set $V$ and identify 
a graph with the set of its edges, then we may regard the set of 
not $i$-connected graphs on $V$ as a simplicial complex.
\\

\begin{definition}
$\Delta_n^i$ is the complex of not 
$i$-connected graphs on $n$ vertices.
\end{definition}

For a graph $G$ and a vertex $v$ we denote by $G-v$ the graph that is
obtained from $G$ by deleting the vertex $v$ from its set of vertices and 
deleting
all edges emerging from $v$ from the set of edges. If $v$ and $w$ are two
distinct vertices of $G$ then we denote by 
$vw$ the two-element set $\{v,w\}$, by $G \setminus vw$ the graph
$(V(G),E(G) \setminus \{vw\})$, and by $G + vw$ the graph 
$(V(G),E(G) \cup \{vw\})$.  Note that (by definition) if $xy \in E(G)$ 
then $G+xy=G$ and if $xy \not\in E(G)$ then $G \setminus xy=G$.
A subset $V' \subseteq V(G)$ of the vertex-set of a
graph $G = (V(G),E(G))$ is called a {\it cutset} if the graph obtained
from $G$ by deleting the vertices in $V'$ and all adjacent edges is
not connected. In particular, a graph is $i$-separable if and only if there is
a cutset of cardinality $i$. A cutset of cardinality $1$ is also called
a {\it cutpoint}.

More generally, one may consider complexes
of not $i$-connected $k$-uniform hypergraphs. Recall that
a {\it $k$-uniform hypergraph} on a vertex set $V$ is a subset $E$ of 
the set of $k$-element subsets $\binom{V}{k}$ of $V$. We will call the
$k$-uniform hypergraphs {\it $k$-graphs} for short. Note that a $2$-graph
is just a graph.  A $k$-graph 
is called {\it $i$-connected} if its underlying $2$-graph is 
$i$-connected.  The underlying $2$-graph of a $k$-graph $E$ is the 
graph on $V$ whose edge set contains a $k$-clique on $\lp v_1,\ldots,v_k \rp$
for each hyperedge $\lp v_1,\ldots,v_k  \rp \in E$.  
\\

\begin{definition}
$\Da^{i}_{n,k}$ is the complex of all not $i$-connected $k$-graphs 
on $n$ vertices.
\end{definition}

Cutsets and cutpoints are defined analogously for $k$-graphs as they
were for graphs.

For the notation related to simplicial complexes and partially ordered sets -- 
posets for short -- used in this paper, we refer the reader to Section 
\ref{notation}.
\\
\\

Let us now review some known results.
For $i=1$ we have that $\Delta_n^1$ and $\Delta_{n,k}^1$ are the complexes of
disconnected graphs, resp., disconnected $k$-graphs. 
The topology of $\Delta_{n,k}^1$ is well understood up to homotopy type.

\begin{proposition} \label{discon} Let $n \geq 2$.  Then
\begin{itemize}
\item[(i)] The complex $\Delta_n^1$ is homotopy equivalent to a wedge 
of $(n-1)!$
spheres of dimension $n-3$. In particular, $\widetilde{H}_{i}
(\Delta_n^1) = 0$ for
$i \neq n-3$ and $\widetilde{H}_{n-3}(\Delta_n^1) \cong \ZZ^{(n-1)!}$.
\item[(ii)] The complex $\Delta_{n,k}^1$ is homotopy equivalent to 
a wedge of spheres of dimensions $n-(k-2)\cdot t -3$, $1 \leq t \leq 
\lfloor \frac{n}{k} \rfloor$.
In particular, the homology of $\Delta_{n,k}^i$ is free and concentrated in
dimensions $n-(k-2)\cdot t -3$, $1 \leq t \leq \lfloor \frac{n}{k} \rfloor$.
\end{itemize}
\end{proposition}

Part (i) follows from well-known properties of partition lattices (see
\cite{Bjo91,BjoWal83,Sta82}) together with the crosscut theorem (see
\cite{Bjo91}). An alternative proof is provided in \cite{Vas93-1}.
Part (ii) was established by Bj\"orner and Welker in \cite{BjoWel92}.
See Theorem 4.5 and Section 7.8 of \cite{BjoWel92} for
exact numerical information on the homology of $\Delta_{n,k}^1$.

The character of the symmetric group for the representation on 
$\widetilde{H}_{n-3}(\Delta_n^1)$ was determined by Stanley in 
\cite{Sta82} in terms of the character of $S_n$ on the homology 
of the partition lattice.  These two characters are equal by an equivariant 
version of the crosscut theorem. The character of the
symmetric group on the homology of $\Delta_{n,k}^1$ was given by
Sundaram \& Wachs \cite{SunWac94}.

Unless otherwise explicitly stated, all homology groups in this paper
have integer coefficients.

\section{Homology and homotopy type of $\Delta_n^2$} \label{homnot2}

The main theorem of this section gives a complete description of the
homotopy type of $\Delta_n^2$.

\begin{theorem} Let $n \geq 3$.
Then $\dee{n}{2}$ has the homotopy type of a wedge of $(n-2)!$ spheres of
dimension $2n-5$.
\label{main}
\end{theorem}

\noindent

\noindent {\sf Remark:} This result was circulated for several months as a 
conjecture.  During that time, the Euler characteristic of $\dee{n}{2}$
was calculated by Rodica Simion \cite{Si}.  
The theorem was proved independently
and simultaneously, almost to the day, by V. Turchin in Moscow, in a
homology version \cite{Vas97} that is equivalent to our result by some 
general arguments from homotopy theory. 

\medskip

For any natural number $k$, let $B_k$ be the Boolean algebra 
on $k$ elements (i.e., the lattice of subsets of a $k$-element set)
and let $\Pi_k$ be the lattice of partitions of a $k$-set into subsets, 
ordered by refinement. It is well-known that $\de(\ov{B_k})$ --- being the
barycentric subdivision of a simplex boundary --- is homeomorphic
to a $(k-2)$-sphere, and that $\de(\ov{\Pi_k}) \simeq \dee{k}{1}$ has 
the homotopy type of a wedge of $(k-1)!$ spheres of dimension $k-3$ 
(see Proposition \ref{discon} (i) and its references).  These
facts imply the following.

\begin{lemma} \label{boolxpart}
$\de(\ov{B_k \times \Pi_k})$ has the homotopy type 
of a wedge of $(k-1)!$ spheres of dimension $2k-3$.
\end{lemma}

\Proof Let $\emptyset$ and $[k]$ be the least element and top element of $B_k$,
and let $1|\cdots |k$ and $|1\cdots k|$ be the least and top elements of
$\Pi_k$. Apply the Homotopy Complementation Formula \ref{homcom} (ii) to 
$p = (\emptyset,|1\cdots k |)$. The set of complements of $p$ in
$B_k \times \Pi_k$ consists of the single element $q = ([k],1| \cdots|k)$. 
Obviously, $\Delta((\hat{0},q)) \cong \Delta(\ov{B_k})$ and 
$\Delta((q, \hat{1})) \cong \Delta(\ov{\Pi_k})$. Then by 
Formula \ref{homcom} (i) we have
$$\Delta(\ov{B_k \times \Pi_k}) \simeq 
\Sigma \big(\Delta(\ov{B_k}) * \Delta(\ov{\Pi_k})\big).$$
Since the join of a wedge of $n$ spheres of dimension $i$ with a wedge of $m$
spheres of dimension $j$ is homotopy equivalent to a wedge of $nm$ spheres of
dimension $i+j+1$ (see for example \cite[Lemma 2.5 (ii)]{BjoWel92}) 
the assertion follows. Recall that suspension can be regarded as a
join with a $0$-sphere and that the join operation is associative.
\qed

Thus, in order to prove Theorem \ref{main} it suffices to demonstrate
that $\dee{n}{2}$ 
is homotopy equivalent to $\de(\ov{B_{n-1} \times \Pi_{n-1}})$.  
In order to state more
precisely what we will prove, we make the following definitions.

\medskip
\begin{definition}
For $x \in \lb n \rb$ and any graph $G$ on $[n]$, $N_G(x)$ is the {\it neighborhood}
of $x$ in $G$, i.e. $N_G(x)=\{y\in [n]: xy\in E(G)\}$,
 and $\pi(x,G)$ is the partition of the set $\lb n \rb \setminus
\lp x \rp$ determined by the connected components of $G-x$.
\end{definition}


\begin{definition}
$\phi:\ov{\Lat(\dee{n}{2})} \rightarrow \ov{B_{n-1} \times \Pi_{n-1}}$ is
the map of posets given by $G \mapsto (N_G(1),\pi(1,G))$, and 
$\phi^\ast:\de(\ov{\Lat(\dee{n}{2})}) \rightarrow \de(\ov{B_{n-1} \times \Pi_{n-1}})$
is the simplicial map induced by $\phi$.
\end{definition}

\noindent
Note that if $G$ is a graph on $[n]$ such that $N_G(1)=\lp 2,\ldots,n \rp$ and
$G-1$ is connected, then $G$ is $2$-connected.  On the other hand, if $N_G(1)=
\emptyset$ and $\pi(1,G)=2|3|\ldots|n$ then $G$ is the empty graph.  Thus 
$\phi$ is well-defined.  It is clear that $\phi$ is order preserving, so 
$\phi^\ast$ is well-defined.  We can now state the key technical result,
from which (in view of Lemma \ref{boolxpart}) 
Theorem \ref{main} follows.  

\begin{lemma} \label{hom:equiv}
The simplicial map $\phi^\ast$ is a homotopy equivalence.
\label{homeq}
\end{lemma}

To prove Lemma \ref{homeq} we use Quillen's Fiber 
Lemma (see Proposition \ref{fibre}).  In our situation this says that
if for each $(S,\pi) \in \ov{B_{n-1} \times \Pi_{n-1}}$ the order complex of
the poset $\phi_{\leq}^{-1}(S,\pi)=\lp G \in \ov{\Lat(\dee{n}{2})}:\phi(G) 
\leq (S,\pi) \rp$ is
contractible, then $\phi^\ast$ is a homotopy equivalence.  
 If $\pi \neq |2 \cdots n|$ then
$\phi^{-1}_{\leq}((S,\pi))$ has a top element, namely the graph $G$ such 
that $1t$ 
is an edge of $G$ for $t \in S$ and $G$ induces the complete graph on each 
block of $\pi$. So assume that $\pi = |2 \cdots n|$.
If $|S| \leq 1$ then there is also a top element in 
$\phi^{-1}_{\leq}((S,\pi))$, namely the graph $G$ which induces a
clique on $\lp 2,\ldots,n \rp$ and has $N_G(1)=S$. If 
$S = \{ 2, \ldots ,n \}$ then $(S, \pi)$
does not lie in the proper part of $B_{n-1} \times \Pi_{n-1}$. 
In summary, it remains to consider the fibers $\phi_{\leq}^{-1}(S,\pi)$
for pairs $(S,\pi)$ such that $\pi = |2 \cdots n|$ and $S \subseteq
\{ 2, \ldots ,n \}$ with $2 \leq |S| \leq n-2$.
To handle these remaining cases, we make the following definitions.

\medskip
\begin{definition}
\begin{itemize}
\item[(1)] For $2 \leq k \leq n-1$, $\de(k)=\lp G \in \dee{n}{2}:N_G(1) \subseteq \lp
2,\ldots,k \rp \rp$.
\item[(2)] For $3 \leq k \leq n-1$, $\de(k-1,k)=\lp G \in \de(k-1):G+1k \in \de(k) 
\rp$.
\end{itemize}
\end{definition}

\noindent
Note that if $(S,\pi)=(\lp 2,\ldots,k \rp,|2 \cdots n|)$ then 
$\de(k)=\phi_{\leq}^{-1}(S,\pi)$.  Also,
$\de(k-1,k)$ consists of those graphs in $\de(k-1)$ which do not become
$2$-connected when the edge $1k$ is added.  

By the above discussion and the
fact that the natural action of $S_n$ on $\Lat(\dee{n}{2})$ is order 
preserving, Lemma \ref{homeq} follows immediately from the next lemma.

\begin{lemma}
For $2 \leq k \leq n-1$, $\de(k)$ is contractible.
\label{dkcon}
\end{lemma}

\noindent
The proof of Lemma \ref{dkcon} proceeds by induction on $k$, 
the case $k=2$ having been handled above.  
The inductive proof is therefore achieved by the combination of the 
following two lemmas.

\begin{lemma}
Let $3 \leq k \leq n-1$.  If $\de(k-1)$ and $\de(k-1,k)$ are contractible, 
then so is $\de(k)$.
\label{redk}
\end{lemma}

\noindent
\begin{Proof}
Let $\star(1k)$ be the subcomplex of $\de(k)$ consisting of graphs that
either contain the edge $1k$ or else can be extended within $\de(k)$
to contain $1k$. Then $\star(1k)$ is a cone with base $\de(k-1,k)$ and
apex $1k$, and we have
$$\de(k) = \de(k-1) \cup \star(1k),$$
$$\de(k-1,k) = \de(k-1) \cap \star(1k).$$
Thus, $\de(k)$ is a union of two contractible complexes with contractible
intersection, and hence $\de(k)$ is itself contractible (see e.g. 
\cite[Lemma 10.3]{Bjo91}).
\end{Proof}

\begin{lemma}
For $3 \leq k \leq n-1$, $\de(k-1,k)$ is contractible.
\label{dkkcon}
\end{lemma}

\smallskip
To prove Lemma \ref{dkkcon} we will use a special case of
Forman's discrete Morse theory (see \cite{fo}, and for this case also \cite{ch}).
The following works for regular cell complexes, but we will only need
the simplicial case.

\smallskip\noindent
\begin{definition}
Let $\Sigma$ be a simplicial complex.

\begin{itemize}
\item[(1)] $D(\Sigma)$ is the digraph whose vertex set is $\Sigma$ and
whose edges are the edges in the Hasse
diagram of $\Lat(\Sigma)\setminus \{\hat{1}\}$, all directed downward.
\item[(2)] For any set $X$ of edges in $D(\Sigma)$, $D_X(\Sigma)$ is the 
digraph
obtained from $D(\Sigma)$ by reversing the direction of the edges in $X$,
so these edges are directed upward while the remaining edges are directed
downward.
\end{itemize}
\end{definition}

Before we can formulate the following lemma we have to recall some basic facts 
about collapsibility (see for example \cite{Bjo91}).
Given a simplicial complex $\Sigma$, a face $\sigma 
\in \Sigma$ is called {\it free} if $\sigma$ is not maximal and is contained 
in a unique maximal face of $\Sigma$. If $\sigma$ is free in $\Sigma$ then 
passing from $\Sigma$ to the
complex $\Sigma \setminus \{ \tau:\tau \supseteq \sigma \}$ is called an 
{\it elementary collapse} of $\Sigma$. If we can obtain a single vertex
by applying a sequence of
elementary  collapses to a complex $\Sigma$,
then $\Sigma$ is called {\it collapsible}. Since it is easily seen that an
elementary collapse of $\Sigma$ is a strong deformation retraction it 
follows that collapsible complexes are contractible. 

\begin{proposition}
Let $\Sigma$ be a simplicial complex.  If $D(\Sigma)$ contains a perfect
matching $M$ such that $D_M(\Sigma)$ is acyclic, then $\Sigma$ is collapsible.
\label{forman}
\end{proposition}

\noindent
\begin{Proof} 
This is a special case of Corollary 3.5 of \cite{fo}, and this case
is easily proved by induction on $|\Sigma|$.  If $\Sigma=\lp \emptyset,\lp x
\rp \rp$ then the claim is clearly true.  If $|\Sigma|>2$, let $x$ be a 
source in $D_M(\Sigma)$, which must exist since $D_M(\Sigma)$ contains no
directed cycle.  It is easy to see that $x$ must be a free face
of $\Sigma$ which is properly contained in a unique face $y \in \Sigma$.  Now
$\Sigma$ is collapsible to the complex obtained by removing $x$ and $y$,
and we can apply the inductive hypothesis.  
\end{Proof}

We call a perfect matching of the type described in Proposition \ref{forman} an 
{\it acyclic perfect matching} on $D(\Sigma)$. Our goal is to
produce an acyclic perfect matching on $D(\de(k-1,k))$.  The following easy
result will be useful.

\begin{lemma}
Let $\Sigma$ be a simplicial complex, let $M$ be a matching on $D(\Sigma)$
and let $F_0 \ra F_1 \ra \ldots \ra F_r \ra F_0$ be a directed cycle in
$D_M(\Sigma)$.  Then there is some dimension $d$ such that $dim(F_i) \in
\lp d,d+1 \rp$ for all $i \in \lb r \rb$.
\label{twolev}
\end{lemma}
\begin{Proof}
If the $F_i$ have more than two distinct dimensions then some $F_i$
must be incident to two upward directed edges.  This contradicts the fact
that $M$ is a matching, and the result follows immediately. 
\end{Proof}

Before proceeding with the proof of Lemma \ref{dkkcon} we make
some technical definitions. 

\begin{definition} Consider separable graphs on the vertex set $[n]$.
\begin{itemize}
\item[(1)] We denote the set of cutpoints of such a graph $G$ by $\Cut(G)$.
\item[(2)] For fixed $k \in \lp 3,\ldots,n-1 \rp$, let
\begin{itemize}
\item[(a)] $I(k):=\lp G \in \de(k-1,k)~|~N_G(1)=\emptyset \rp$.
\item[(b)] $J(k):=\lp G \in \de(k-1,k)~|~N_G(1) \neq \emptyset \mbox{ and } 
\Cut(G+1k) \neq \lp 1 \rp \rp$.
\item[(c)] $F(k):=\lp G \in \de(k-1,k)~|~\Cut(G+1k)=\lp 1 \rp \rp$.
\end{itemize}
\end{itemize}
\end{definition}

\noindent Note that $\de(k-1,k)$ is the disjoint union of 
$I(k)$, $J(k)$ and $F(k)$, 
and that both $I(k)$ and $I(k) \cup J(k)$ are subcomplexes of $\de(k-1,k)$.

\smallskip

The following lemma implies Lemma \ref{dkkcon}, and
therefore completes the proof of Theorem \ref{main}.

\begin{lemma}
For any $k \in \lp 3,\ldots,n-1 \rp$, $D(\de(k-1,k))$ admits an acyclic
perfect matching.
\label{apm}
\end{lemma}
\noindent
\begin{Proof} This proof will be carried out in three steps. We will 
construct an acyclic perfect matching first for $D(I(k))$,
then for $D(I(k) \cup J(k))$, and finally for $D(\de(k-1,k))$.

\medskip

\noindent {\sf Step 1:} $D(I(k))$ admits an acyclic perfect matching. 

Note 
that $I(k)$ contains a unique maximal face, namely the complete graph on 
$\lp 2,\ldots,n \rp$. Thus $I(k)$ is a simplex and it is easy to see that 
the matching $M=\lp G+23 \rightarrow G \setminus 23 | G \in I(k) \rp$ is 
an acyclic perfect matching on $D(I(k))$.

\medskip

\noindent {\sf Step 2:} $D(I(k) \cup J(k))$ admits an acyclic perfect matching.

It suffices
to show that there exists a matching $M^\ast$ consisting of edges between
elements of $J(k)$ which covers all the elements of $J(k)$, and such that 
$D_{M^\ast}(I(k) \cup J(k))$ is acyclic. If $M^\ast$ is such a 
matching, let $M^\circ$ be an acyclic perfect matching on $D(I(k))$ and set
$M=M^\ast \cup M^\circ$.  Then $M$ is a perfect matching on 
$D(I(k) \cup J(k))$ which contains no edges between $I(k)$ and $J(k)$, so 
that any directed cycle in $D_M(I(k))$ cannot cover points from both $I(k)$
and $J(k)$.  It follows immediately that $M$ is acyclic.

\smallskip

Now let $G \in J(k)$ and let $c \in \Cut(G+1k)$, $c\neq 1$.  
Let $x=\min\lp N_G(1) \rp$.  
If $xk  \in E(G)$ then clearly $G\setminus xk \in J(k)$.  If 
$c \not\in \lp x,k \rp$ then since $1k$ and
$1x$ are edges of $G+1k$, $x$ and $k$ lie in the same connected component
of $(G+1k)-c$.  If $c \in \lp x,k \rp$ then clearly $c$ is a cutpoint of 
$G+xk+1k$.  In any case, $c \in \Cut(G+xk+1k)$ and $G+xk \in J(k)$.

\smallskip

Let $M^\ast$ consist of all edges $G+xk \rightarrow G \setminus xk$, where $x$ is
determined as above.  Clearly $M^\ast$ is a matching which covers all points
in $J(k)$.  Assume for contradiction that $A_1 \ra B_1 \ra A_2 \ra B_2 \ra
\ldots \ra B_r \ra A_1$ is a directed cycle in $D_{M^\ast}(I(k) \cup J(k))$.
Clearly all the $A_i$ and all the $B_i$ are in $J(k)$, and by Lemma
\ref{twolev} we may assume that for each $i$ there are edges $\alpha_i$ and 
$\beta_i$ such that $B_i=A_i+\alpha_i$ and $A_{i+1}=B_i \setminus \beta_i$.  
Thus
$A=A+\alpha_1 \setminus \beta_1+\ldots +\alpha_r\setminus \beta_r$ and 
$\lp \alpha_i \rp=\lp \beta_i \rp$.
By the definition of $M^\ast$, no $\alpha_i =x_i k$ contains $1$, so no $\beta_i$
contains $1$.  It follows that $N_{A_1}(1)=N_{B_1}(1)=N_{A_2}(1)=\ldots=
N_{B_r}(1)$. By the choice of the $x_i$'s 
this forces $\alpha_1=\alpha_2=\ldots=\alpha_r$, which is
clearly impossible.

\medskip

\noindent {\sf Step 3:} $D(\de(k-1,k))$ admits an acyclic perfect matching.

As in Step 2,
it suffices to produce a matching $M^\ast$ on edges connecting elements of
$F(k)$ which covers all points in $F(k)$ and such that $D_{M^\ast}(\de(k-1,k))$
is acyclic.

\smallskip

Let $G \in F(k)$.  Then $(G+1k)-1$ splits into connected components $C_1,
\ldots,C_s$ such that for each $i \in \lb s \rb$ the subgraph of $G+1k$
induced on $V(C_i) \cup \lp 1 \rp$ is $2$-connected.  We may assume that
$n \in V(C_1)$.  Note that since $k<n$, $1n \not\in G+1k$.  Define $S(G)$ to 
be the set of all $x \in V(C_1) \cap N_{G+1k}(1)$ such that there is a path
$P=1,x,\ldots,n$ in $G+1k$ with $P \cap N_{G+1k}(1)=\lp x \rp$.

\smallskip

We claim that $|S(G)|>1$.  Indeed, let $1,x,\ldots,n$ be a shortest path from
$1$ to $n$ in $G+1k$.  Clearly $x \in S(G)$.  Since the subgraph of $G+1k$
induced on $V(C_1) \cup \lp 1 \rp$ is $2$-connected and $x \neq n$, there 
exists a path from $1$ to $n$ in this graph which does not contain $x$.  
Let $1,y,\ldots,n$ be a shortest such path.  Then $y \in S(G)$.

\smallskip

Let $x,y$ be the two smallest elements of $S(G)$.  If $xy \not\in G$ then
clearly $G+xy \in F(k)$ and $S(G+xy)=S(G)$.  Now assume $xy \in G$ and let
$H$ be the subgraph of $G \setminus xy+1k$ induced on $V(C_1) \cup \lp 1 \rp$.  If 
$d$ is a
cutpoint of $H$ then $x$ and $y$ are in different components of $H-d$ 
(otherwise $d$ is a cutpoint of the subgraph of $G+1k$ induced on $V(C_1)
\cup \lp 1 \rp$).  However, there is a cycle $1,x,\ldots,y,1$ in $H$.  Thus
there is no such cutpoint $d$ and $H$ is $2$-connected.  It follows that
$G \setminus xy \in F(k)$ and $S(G \setminus xy)=S(G)$.

\smallskip

Now, let $M^\ast$ consist of the edges $G+xy \ra G \setminus xy$ where $x,y$
are determined as above.  Then $M^\ast$ is a matching which consists of
edges connecting points in $F(k)$ and covers all points in $F(k)$.  It
remains to show that $D_{M^\ast}(\de(k-1,k))$ is acyclic.

\smallskip

Assume for contradiction that $A_1 \ra B_1 \ra A_2 \ra \ldots \ra B_r \ra A_1$
is a directed cycle in $D_{M^\ast}(\de(k-1,k))$.  As in Step 2, we may assume
that there are edges $\alpha_i$ and $\beta_i$ such that $B_i=A_i+\alpha_i$,
$A_{i+1}=B_i \setminus \beta_i$ and $\lp \alpha_i \rp=\lp \beta_i \rp$.

\smallskip

By the definition of $M^\ast$, each $\alpha_i$ connects two elements of
$N_{B_i+1k}(1)=N_{A_i+1k}(1)$, so no $\beta_i$ contains $1$.  Thus $N_{A_i}(1)
=N_{B_j}(1)$ for all $i,j$, and each $\beta_i$ connects two elements of
$N_{B_i+1k}(1)=N_{A_{i+1}+1k}(1)$.  Write $\alpha_1=xy$.  Then $\beta_1 \neq
xy$, and in $A_2+1k$, $x$ and $y$ are still the two smallest neighbors of $1$
which are contained in paths from $1$ to $n$ which intersect $N_{A_2+1k}(1)$
exactly once.  Thus $\alpha_2=xy=\alpha_1$, giving the desired
contradiction.  
\end{Proof}

\section{The character for the action of $S_n$ on 
$\ti{H}_{2n-5}(\dee{n}{2})$}
\label{rep} \label{homrep}

In view of Theorem \ref{main} it is natural to investigate the
representation of the symmetric group $S_n$ on the only non-zero homology
group of $\dee{n}{2}$, induced by the obvious action. In this section we
consider homology with complex coefficients, hence all representations
are over $\CC$.  In many of the computations below, we actually determine 
character values for the representation of $S_n$
on the only non-zero homology group of $\Delta(\ov{\Lat(\dee{n}{2})})$, which
is easily seen to be the same as the representation described above.

\medskip
\begin{definition}
\begin{itemize}
\item[(i)] We denote by $\omega_n^2$ the character of $S_n$ given by $g \mapsto 
Trace(g,\ti{H}_{2n-5}(\dee{n}{2}))$.
\item[(ii)] Let $C_n$ be a cyclic subgroup of $S_n$ generated by a full
$n$-cycle. We denote by $\lie_n$ the character of $S_n$ induced from
the character on $C_n$ which takes the value $e^{\frac{2 \pi i}{n}}$ on a
fixed generator. It is well known (see e.g. \cite[Chapter 8]{Reu93}) that 
$\lie_n$ is the
character of $S_n$ on the multigraded piece of the free Lie algebra generated
by $n$ variables.
\end{itemize}
\end{definition}

For the rest of this section we let $S_{n-1}$ be the stabilizer of the point
$1$ in the natural action of $S_n$ on the set $[n]$.

\begin{theorem} \label{char}
The character $\omega_n^2$ is given by 
$$\omega_n^2 = \lie_{n-1} \uparrow_{S_{n-1}}^{S_n}
- \lie_{n}.$$ 
\end{theorem} 
The proof will follow a sequence of lemmas establishing the main steps.

\begin{lemma} \label{first}
If $g \in S_{n-1}$ then $\omega_n^2(g)=\lie_{n-1}(g).$ 
\label{resg1}
\end{lemma}
\begin{Proof} 
It is easily seen that the map $\phi:\ov{\Lat(\dee{n}{2})} \rightarrow
\ov{B_{n-1} \times \Pi_{n-1}}$,
defined in the previous section, commutes with the actions of $S_{n-1}$ on the
two posets. Thus the induced map on homology is $S_{n-1}$-equivariant and
is an $S_{n-1}$-module isomorphism by Lemma \ref{hom:equiv}. Thus,
the characters of $S_{n-1}$ on the homology of $\Delta_n^2$ and on the homology
of $\Delta(\ov{B_{n-1} \times \Pi_{n-1}})$ coincide. By an equivariant
version of Proposition \ref{homcom} (see \cite{Wel90-1}),
$\Delta(\ov{B_{n-1} \times \Pi_{n-1}})$ has the $S_{n-1}$-homotopy type of 
$\Sigma (
\Delta(\ov{B_{n-1}}) * \Delta(\ov{\Pi_{n-1}})),$ where the group $S_{n-1}$ acts
diagonally on $\Delta(\ov{B_{n-1}}) * \Delta(\ov{\Pi_{n-1}})$. 
Thus the character of 
$S_{n-1}$ on the homology of $\Delta(\ov{B_{n-1} \times \Pi_{n-1}})$ is given by
the product of the characters of $S_{n-1}$ on 
$\widetilde{H}_*(\Delta(\ov{B_{n-1}}))$ and 
$\widetilde{H}_*(\Delta(\ov{\Pi_{n-1}}))$. The character of $S_{n-1}$ on
$\widetilde{H}_*(\Delta(\ov{B_{n-1}}))$ is rather easily seen to be the
sign-character of $S_{n-1}$ (see \cite{Sta82}). The character of
$S_{n-1}$ on $\widetilde{H}_*(\Delta(\ov{\Pi_{n-1}}))$ was
determined in \cite{Sta82} as 
$sign_{n-1} \cdot \lie_{n-1}$. This implies the assertion. 
\end{Proof}

\noindent
Since every element of $S_n$ which has a fixed point is conjugate to an
element of $S_{n-1}$, it remains to determine $\omega_n^2(g)$ for all 
fixed-point-free $g \in S_n$.

\medskip
\begin{definition}
Let $g \in S_n$.  We denote by $L^g$ the poset of faces of 
$\dee{n}{2}$ which are
fixed by $g$, and by $g^\ast$ the element of $S_{n+1}$ which
fixes $n+1$ and acts as $g$ does on $\lb n \rb$.
\end{definition}

\noindent
Write $\ho$ for the empty graph in $L^g$, which is the unique minimum 
element of $L^g$, and for any poset $P$ let $\mu_{P}$ be the M\"obius 
function on $P$.

\begin{lemma} 
For $g \in S_n$, $\omega_n^2(g)=\displaystyle{\sum_{G \in L^g}\mu_{L^g}(\ho,G)}$.
\label{hall}
\end{lemma}

\noindent
\begin{Proof} It is well-known (see e.g. \cite[(13.5)]{Bjo91})
that if a group
acts on a bounded poset $P$ then for any group element $g$ we have
\[
\mb{P^g}=\sum_i(-1)^iTr(g,\ti{H}_i(\de(\ov{P}))).
\]
In the case under consideration, the only nonzero reduced homology group
is the one in dimension $2n-5$, so the lemma follows immediately from the
definition of the M\"obius function.  
\end{Proof}

The next two lemmas will be used to determine $\omega_n^2(g)$ when $g$ is 
fixed-point-free.  

\begin{lemma}
Let $G$ be a graph whose automorphism group acts transitively on $V(G)$.  
If $G$ is connected then $G$ is $2$-connected.
\label{cyctrans}
\end{lemma}

\begin{Proof}
Let $v$ be a leaf of some spanning tree in the connected graph $G$.
Then $v$ is not a cutpoint. Since Aut($G$) is transitive on vertices there 
cannot be any other cutpoints. Hence $G$ is $2$-connected. 
\end{Proof} 

\begin{lemma}
Let $g \in S_n$ be fixed-point-free.  Write $g$ as a product of disjoint 
cycles, $g=g_1 \ldots g_r$.  Let $V_i=supp(g_i)$.  Let $G \in L^g$ be 
connected and let $x \in \Cut(G)$ with $x \in V_j$.  Then there exists some
connected component $C$ of $G-x$ such that $V_j \setminus \lp x \rp 
\subseteq C$ and $C \cap V_i \neq \emptyset$ for all $i \in \lb r \rb$.
\label{bigcomp}
\end{lemma}

\noindent
\begin{Proof} 
Let $G_j$ be the graph on $V_j$ such that an edge $yz$ is in $E(G_j)$
if $yz \in E(G)$ or if there is a path $P$ from $y$ to $z$ in $E(G)$ such that
$P \cap V_j=\lp y,z \rp$.  Since $G$ is connected, so is $G_j$.  Also, the
group generated by $g_j$ is a group of automorphisms of $G_j$ which acts 
transitively on $V_j$.  By Lemma \ref{cyctrans},
$G_j$ is $2$-connected.  It follows that all elements of $V_j \setminus \lp x
\rp$ are in the same connected component of $G-x$.  Now for $i \neq j$, let
$P$ be a path of shortest length connecting some $y \in V_i$ with some $z
\in V_j$.  If $z=x$ replace $P$ with $g(P)$.  Now $P$ contains no vertices
from $V_i \cup V_j$ other than $y$ and $z \neq x$.  Thus $P$ is a path in
$G-x$ and $y$ lies in the component of $G-x$ containing $V_j \setminus \lp 
x \rp$.  \end{Proof}

We can now determine the values of $\omega_n^2$ on fixed-point-free 
elements of $S_n$.

\begin{lemma} \label{second}
Let $g \in S_n$ be fixed-point-free.  
Then $\omega_n^2(g)=-\omega_{n+1}^2(g^\ast)$.
\label{wfpf}
\end{lemma}
\begin{Proof} 
As usual we write $\ho$ for the empty graph.  By Lemma \ref{hall} we have
\[
\omega_{n+1}^2(g^\ast)=\sum_{G \in L^{g^\ast}}\mega.
\]
Let $M^g$ be the poset of all graphs on $\lb n \rb$ which are fixed by $g$.
Note that if $G \in L^{g^\ast}$ then $G - (n+1) \in M^g$.  For $F \in M^g$
let $D(F)$ be the set of all $G \in L^{g^\ast}$ such that $G-(n+1)=F$.  We
have
\[
\omega_{n+1}^2(g^\ast)=\sum_{F \in M^g}\sum_{G \in D(F)}\mega.
\]
Any $G \in L^{g^\ast}$ is a union of $\langle g^\ast \rangle$-orbits on
${{\lb n+1 \rb} \choose {2}}$.  Let $o(G)$ be the number of such orbits.
It is easy to see that $\mega=(-1)^{o(G)}$.  Let $p(G)$ be the number of
such orbits containing edges covering the point $n+1$.  Applying the
previous argument to $M^g$, we get for any $F \in M^g$
\[
\sum_{G \in D(F)}\mega=\meg\sum_{G \in D(F)}(-1)^{p(G)}.
\]

We will examine this sum for each $F \in M^g$, looking separately at the
cases where $F$ is disconnected, connected but not $2$-connected, and
$2$-connected.  Write $g$ as a product of disjoint cycles, $g=g_1 \ldots g_r$
and let $V_i=supp(g_i)$.  Note that if $v \in V_i$ and $G \in L^{g^\ast}$
with $\lp v,n+1 \rp \in G$, then the $\langle g^\ast \rangle$-orbit 
containing $\{v,n+1\}$ consists of the edges $\{w,n+1\}$ for all 
$w \in V_i$, and is contained in $E(G)$.  Also, $p(G)$ is simply the number 
of such orbits.  Let $O(g)$ be the set of all such orbits, and for $S
\subseteq O(g)$ let $G(S)$ be the graph induced on the edges which are
contained in elements of $S$.  For 
$F \in M^g$ define
\[
\Sigma(F):=\lp S \subseteq O(g): F \cup G(S) \in L^{g^\ast} \rp.
\]
Note that $\Sigma(F)$ is a simplicial complex on $O(g)$.  Let 
$P(F)=\Lat(\Sigma(F)) \setminus \lp \hat{1} \rp$.  By the above arguments 
we have
\[
\sum_{G \in D(F)}\mega=\meg\sum_{S \in P(F)}\mu_{P(F)}(\hat{0},S).
\]
\medskip

We now examine the three cases.

\medskip

\noindent {\sc Case 1:} $F$ is not connected.  

Then $P(F)$ is the Boolean algebra on
$O(g)$, since $n+1$ is a cutpoint of $F \cup G(S)$ for all $S \subseteq O(g)$.
It follows immediately that
\[
\sum_{G \in D(F)}\mega=0.
\]

\medskip

\noindent {\sf Case 2:} $F$ is connected but not $2$-connected.  

We will use the block decomposition described in Proposition \ref{sep}
of the following section.  Given a connected but 
not $2$-connected graph $F \in M^g$, let $T(F)$ be the bipartite graph whose
vertices are the vertices of $F$ and the blocks of $F$, with $\{ v,W_i\}$ 
an edge if and only if $v \in W_i$.  It is easy to see that $T(F)$ is a tree
and that $\langle g \rangle$ is a group of automorphisms of $T(F)$ which
preserves each part of the given bipartition.  It follows that $g$ fixes
a vertex of $T(F)$ (see \cite{Lov93}).  Since $g$ fixes no vertex of $F$, $g$
must fix some block $B$ of $F$.  This means that there is some nonempty
$J \subset \lb r \rb$ such that $B=\cup_{j \in J}V_j$.  Let $S$ be the set
of all orbits in $O(g)$ which contain edges that include vertices in $B$.  
We will show that every maximal element of $P(F)$
contains $S$, from which it follows immediately that
\[
\sum_{G \in D(F)}\mega=0.
\]

Let $G \in D(F)$ and let $c \in \Cut(G)$.  Since $F$ is connected, $c \neq 
n+1$.  Also, if $N_G(n+1) \neq \emptyset$ then since $g$ is fixed-point-free
$c$ must be a cutpoint of $F$.  If $N_G(n+1)=\emptyset$ then every $x \in
\lb n \rb$ cuts $G$, so in any case we may assume $c \in \Cut(F)$.  Let $c
\in V_i$.  By Lemma \ref{bigcomp}, there is some connected component $C$
of $F-c$ which contains $V_i \setminus \lp c \rp$ and at least one element
of each $V_j$.  Since $B$ is $2$-connected and $B \cap C \neq \emptyset$,
we must have $B \subseteq C$.  

We will now show that $c$ must be a cutpoint
of $G \cup S$.  If $N_G(n+1)=\emptyset$ then adding $S$ to $G$ simply moves
the previously isolated point $n+1$ into the connected component of $G-c$
which contains $C$.  However, there is a component of $F-c$ besides $C$, 
which remains separated from $C$ in $G-c$.  Now, assume that $N_G(n+1) \neq
\emptyset$.  Then there exists some set $I$ such that 
$N_G(n+1)=\cup_{i \in I}V_i$.  The component of $G-c$ containing $C$
contains elements of each $V_i$, and it follows that $n+1$ must also be in
this component.  Thus, adding $S$ to $G$ does not reduce the number of 
components of $G-c$.

\medskip

\noindent {\sf Case 3:} $F$ is $2$-connected.  

In this case the only $G \in D(F)$ is that
for which $N_G(n+1)=\emptyset$.  Indeed, since each $V_i$ has at least two
elements we cannot have $|N_G(n+1)|=1$, and the claim follows.  Thus
\[
\sum_{G \in D(F)}\mega=\mu_{M^g}(\ho,F).
\]

\smallskip
Let $K^g$ be the set of $2$-connected graphs in $M^g$.  
Combining the information from the three cases we have shown that
\[
\omega_{n+1}^2(g^\ast)=\sum_{F \in K^g}\mu_{M^g}(\ho,F).
\]
By definition of the M\"obius function and the fact that $M^g$ has a
maximum element, we have
\[
\sum_{F \in K^g}\mu_{M^g}(\ho,F)=-\sum_{F \in L^g}\mu_{L^g}(\ho,F)
=-\omega_n^2(g),
\]
and the proof is complete. 
\end{Proof}

\medskip

\noindent
{\sf Proof of Theorem \ref{char}:}
Set $\rho_n = \lie_{n-1} \uparrow_{S_{n-1}}^{S_n} - \lie_n$.
We must show that $\rho_n(g) = \omega_n^2(g)$ for all $g \in S_n$. 

By the definition of induced characters, if $g \in S_n$ is 
not the product of disjoint cycles of
the same length then $\lie_n(g) = 0$. We will assume from now on that 
any $g \in S_n$ which fixes a point is contained in $S_{n-1}$
(so by our convention it fixes the point $1$). 

By the definition of induced characters and Theorem \ref{main}, we have
$$\rho_n(\id) = (n-2)!~ [S_n:S_{n-1}] - (n-1)! = (n-2)! = \omega_n^2(\id).$$  

If $g \neq \id$ and $g$ has at least two fixed points, then $\lie_{n-1}(g)
= \lie_n(g) = 0$, so $\rho_n(g) = \omega_n^2(g)$ by Lemma \ref{first}. 

If $g \neq \id$ has exactly one fixed point, then $\lie_n(g) = 0$.
For $h \in S_n$ we have $g^h := h^{-1}gh \in S_{n-1}$ if and only if $h \in S_{n-1}$.
By the definition of induced characters and Lemma \ref{first},
$$\rho_n(g) = \lie_{n-1} \uparrow_{S_{n-1}}^{S_n}(g) =
\frac{1}{(n-1)!} \sum_{h \in S_{n-1}} \lie_{n-1}(g^h) = \lie_{n-1}(g) = 
\omega_n(g).$$

If $g \in S_n$ has no fixed points then $\lie_{n-1} \uparrow_{S_{n-1}}^{S_n} 
(g) = 0$ and $\rho_n(g) = - \lie_n(g)$. As before, let $g^*$ be the element of $S_{n+1}$
which fixes $n+1$ and acts as $g$ does on $[n]$. We have shown above that
$\omega_{n+1}^2(g^*) = \lie_n(g)$. Hence, by Lemma \ref{second},
$\rho_n(G) = \omega_n^2(g)$.
\qed

\medskip

According to Vassiliev, the number of linearly independent knot
invariants of bi-order $(n,n-1)$, modulo lower bi-order invariants,
is bounded from above by the multiplicity of the
trivial representation in the restriction of $\omega_n^2$ to the cyclic group
$C_n$ generated by $(12 \cdots n)$. See \cite{Vas97} for all details.
As a corollary of Theorem \ref{char} we obtain a formula for
this multiplicity.
We write $\langle \xi,1 \rangle$ for the 
multiplicity of the trivial character in any character $\xi$ of $C_n$. 

\begin{corollary} \label{trivmult}
\[
\langle \omega_n^2 \downarrow^{S_n}_{C_n},1 \rangle=(n-2)! - 
\frac{1}{n}\sum_{d|n}\mu(d)\phi(d) (\frac{n}{d}-1)!
\hskip2pt d^{\frac{n}{d}-1}
\]
\end{corollary}
\begin{Proof}
As a consequence of a result by Hanlon \cite{Han82} (see also \cite{Sta82}), it is straightforward 
to show that
\[
\langle \lie_n\downarrow^{S_n}_{C_n},1 \rangle=
\frac{1}{n}\sum_{d|n}\mu(d)\phi(d) (\frac{n}{d}-1)!
\hskip1pt d^{\frac{n}{d}-1},
\]
where $\mu$ is the usual number-theoretic M\"obius function and $\phi$ is 
Euler's function. On the other hand,
$$\langle \lie_{n-1}\uparrow_{S_{n-1}}^{S_n}
\downarrow_{C_n}^{S_n},1 \rangle = \frac{1}{n} \sum_{g \in C_n}
\lie_{n-1} \uparrow_{S_{n-1}}^{S_n}(g) = \frac{1}{n} \hskip2pt \lie_{n-1}
\uparrow_{S_{n-1}}^{S_n} (\id) =  (n-2)!.$$ Now the assertion follows 
immediately from Theorem \ref{char}. 
\end{Proof}

The values of $w_n = \langle \omega_n^2\downarrow^{S_n}_{C_n},1 \rangle$ for small $n$
are given in the table below.

\medskip

\begin{center}
$\begin{array}{|c||c|c|c|c|c|c|c|c|c|}
\hline
n & 3 & 4 & 5 & 6 & 7 & 8 & 9 & 10 & 11 \\
\hline \hline
w_n & 1 & 1 & 2 & 6 & 18 & 96 & 564 & 4,072 & 32,990\\
\hline
\end{array}$
\end{center}

\centerline{{\sf Table 1:} Multiplicity $w_n$ of the trivial character in $\omega_n^2
\downarrow_{C_n}^{S_n}$}

\medskip

The character $\omega_n^2$ and the tensor product of $\omega_n^2$ with the 
sign character have recently appeared in various different
settings, see Section \ref{questionrep}.

\section{The lattice of block-closed graphs}
\label{closed}

In this section we will obtain information on the topology of 
$\Delta_{n,k}^2$ by producing a lattice $\snk$ such that $\Delta(\ov{\snk})$
is homotopy equivalent to $\Delta_{n,k}^2$ and examining the structure of
$\snk$. For lattice and poset terminology not explained in Section \ref{notation}
we refer to \cite{Sta86}.

We begin by recalling some elements of the well known structure theory of 
separable graphs, which appears e.g. in \cite{Lov93}.

\begin{definition}
Let $G$ be any graph.  A {\em block} of $G$ is a subset $W$ of $V(G)$ such
that the subgraph of $G$ induced on $W$ is $2$-connected or $W$ is a 
singleton or a pair of points connected by an edge, and the subgraph
of $G$ induced on any proper superset of $W$ is separable.
We will say that $G$ is {\em block-closed} if the subgraph induced
on each block is a clique.
\end{definition}

Given a graph $G$, say that $e\equiv e'$ for two of its edges $e$ and
$e'$ if they both lie in some circuit of $G$. This is easily seen to
be an equivalence relation on $E(G)$. If $W$ is the set of nodes 
underlying an equivalence class then $W$ is a block, and all non-singleton blocks
correspond to equivalence classes of edges in this way. {From} this it is
easy to derive the following basic facts about the ``block decomposition''
of $G$, see \cite{Lov93} for more details.

\begin{proposition}
Let $G$ be a graph.  Then there exists a unique decomposition
of $V(G)$ into blocks $W_1,\ldots,W_r$, and if $i \neq j$ we have
$|W_i \cap W_j| \leq 1$.  Moreover, if $B_G$ is the graph with vertex
set $\lp w_1,\ldots,w_r \rp$ such that $\lp w_i,w_j \rp \in E(B_G)$ if
and only if $|W_i \cap W_j|=1$, then $B_G$ is a forest (that is, $B_G$
contains no cycles).
\label{sep}
\end{proposition}

Note that if $K$ is a $k$-graph with underlying graph $G$,
then every block of $G$ has size at least $k$ or is a single vertex.

\begin{definition}
Let $K$ be a $k$-graph with underlying graph $G$, and let
$W_1,\ldots,W_r$ be the blocks of $G$.  We define $K^*$ to be the
$k$-graph which induces the complete $k$-graph on each $W_i$ and
contains no other hyperedges.  We also define $\snk$ to be the poset
of all graphs on vertex set $\lb n \rb$ in which every block is either
an isolated vertex or a clique of size at least $k$, ordered by inclusion.
\end{definition}

The first part of the following lemma is immediate from the definition,
and the second follows via a standard argument for closure operators
on lattices.

\begin{lemma}
\begin{itemize}
\item[(i)] The map $K \mapsto K^*$ defines a closure operator
on ${\Lat(\Delta_{n,k}^2)}$ whose image is isomorphic to $\snk$.
\item[(ii)] $\snk$ is a lattice.
\label{clop}
\end{itemize}
\end{lemma}

The meet operation in the lattice $\snk$ is intersection of edge-sets
followed by deletion of the edges in all blocks of size smaller than $k$.
Note that the elements of $\Sigma_{n,2}$
are the block-closed graphs, and that
we have a tower of embeddings as subposets (not sublattices):
$$\snk\subseteq \cdots \subseteq \Sigma_{n,3} \subseteq \Sigma_{n,2}.$$
Hence, in view of the following result the topology of all the complexes
$\Delta_{n,k}^2$ is encoded into the lattice $\Sigma_{n,2}$ of block-closed graphs.

\begin{theorem}
The complexes $\Delta_{n,k}^2$ and $\Delta(\ov{\snk})$ are homotopy equivalent.
\label{dnksnk}
\end{theorem}

\begin{Proof} $K^*$ is the complete graph (the top element of $\snk$)
if and only if $K$ is $2$-connected. Hence,
the map $K \mapsto K^*$ restricts to a closure operator
on $\ov{\Lat(\Delta_{n,k}^2)}$ whose image is isomorphic to $\osnk$.
The theorem then follows from Corollary \ref{closure}.
\end{Proof}

We will now investigate the structure of $\snk$.  The next two lemmas follow
immediately from the definition of $\snk$.  We write $\hat{0}$ for the
empty graph, which is the minimum element of $\snk$, and $\hat{1}$ for the
complete graph, which is its maximum.

\begin{lemma}
Let $G,H \in \snk$.  Then $G$ covers $H$ if and only if one of the following
conditions holds:
\begin{itemize}
\item[(i)] $E(G) \setminus E(H)$ is a clique on $k$ vertices belonging to $k$
pairwise different components of $H$.
\item[(ii)] $E(G) \setminus E(H)$ is a complete bipartite graph on parts
$A$ and $B$, and there is a vertex $v$ such that $A \cup \lp v \rp$ and
$B \cup \lp v \rp$ are blocks in $H$.
\item[(iii)] Only if $k>2$: $E(G) \setminus E(H)$ is a star (that is, a connected graph
with at most one vertex of degree more than one), and the vertices
of degree one in this star form a block in $H$ belonging to a component
of $H$ distinct from that of the center of the star.
\end{itemize}
\label{covsnk}
\end{lemma}

The three types of coverings can informally be described as follows:
\begin{itemize}
\item[(i)] select a vertex from each of $k$ pairwise disjoint components of $H$
and then create a $k$-clique on these vertices;
\item[(ii)] complete the union of two overlapping blocks of $H$ to a clique;
\item[(iii)] for $k>2$: select a block and a vertex from different components of $H$
and complete their union to a clique.
\end{itemize}

The lattices $\snk$ are neither upper nor lower semimodular. However, they
exhibit a recursive structure on lower intervals, and certain upper intervals
are upper semimodular, as the following lemma shows.

\begin{lemma} Let $G\in \snk$.
\begin{itemize}
\item[(i)] If $G$ has $r$ non-singleton blocks of sizes $m_1, \dots , m_r$ then
the interval $[\hat{0},G]$ is isomorphic to the direct product
$\Sigma_{m_1,k} \times \cdots \times \Sigma_{m_r,k}$
\item[(ii)] If $G$ is connected then the interval $[G,\hat{1}]$ is isomorphic
to a direct product of partition lattices. More precisely, suppose that $G$
has $s$ cutpoints and that the $i$-th cutpoint lies in $t_i \geq 2$ blocks.
Then, $[G,\hat{1}] \cong \Pi_{t_1} \times \cdots \times \Pi_{t_s}$.
\end{itemize}
\label{interv}
\end{lemma}

The following description of the coatoms of $\snk$, that is, the elements 
which are covered by $\hat{1}$, follows immediately from the two
preceding lemmas.

\begin{lemma}
Let $M$ be a coatom of $\snk$. Then one of the following conditions holds:
\begin{itemize}
\item[(i)] $M$ is connected and has two blocks of size $l,m$ with $k \leq
l \leq m \leq n-k+1$ and $l+m=n+1$.  In this case, the interval 
$\lb \hat{0},M \rb$ is isomorphic to $\Sigma_{l,k} \times \Sigma_{m,k}$.
\item[(ii)] $M$ consists of an $(n-1)$-clique and an isolated vertex.  In
this case, $k>2$ and the interval $\lb \hat{0},M \rb$ is isomorphic to
$\Sigma_{n-1,k}$.
\end{itemize}
\label{cosnk}
\end{lemma}

For any graph $G$ let $c(G)$ be the number of connected components, and 
$b(G)$ the number of blocks of size $\geq 2$.

\begin{theorem}
\begin{itemize}
\item[(i)] The lattice $\Sigma_{n,2}$ is graded with rank function
$$\rho(G)=2n-2c(G)-b(G).$$
In particular, its length is $\rho(\hat{1})=2n-3$.
\item[(ii)] The lattice $\Sigma_{n,3}$ is graded with rank function
$$\rho(G)=n-c(G)-b(G).$$
In particular, its length is $\rho(\hat{1})=n-2$.
\item[(iii)] If $k>3$ and $n<2k-1$, then $\snk$ is isomorphic to the 
lower-truncated Boolean algebra 
$\{A\subseteq [n]: |A|\geq k\}\cup \{ \emptyset \}.$
In particular, $\snk$ is graded of length $n-k+1$.
\item[(iv)] If $k>3$ and $n \geq 2k-1$, then $\ell$ is the length of a maximal
chain of $\snk$ if and only if $$\ell =(n-2)-t(k-3), \mbox{ for some }
1\leq t\leq \lfloor \frac{n-1}{k-1} \rfloor .$$
In particular, $\snk$ is of length $n-k+1$ and is not graded.
\item[(v)] If $k>3$ then $G \in \ov{\snk}$ is contained in a chain of
length $n-k+1$ if and only if $G$ consists of a clique of size $l \geq k$
and $n-l$ isolated vertices.
\end{itemize}
\label{chsnk}
\end{theorem}

\begin{Proof}
For claims (i) and (ii) it suffices to check that the given rank functions
increase by $1$ for each type of covering given in Lemma \ref{covsnk}
and take value zero at the empty graph. Claim (iii) is clear from
the definition.

Claims (iv) and (v) are implied by the following description of the
maximal chains in $\snk$. We will here view $\snk$ as a subposet of
$\Sigma_{n,3}$, and we let $\rho$ denote the restriction of the rank 
function of claim (ii) from $\Sigma_{n,3}$ to $\snk$.

A maximal chain from $\hat{0}$ to $\hat{1}$ in $\snk$ is a sequence of
covering steps. By Lemma \ref{covsnk} there are three possibilities
for each step. The rank function $\rho$ will increase by $1$ for
coverings of types (ii) or (iii), and by $k-2$ for coverings of type (i).
Hence, the length of a maximal chain must be $n-2-t(k-3)$, where
$t$ is the number of covering steps of type (i). Note that $t\geq 1$
since the first covering in the chain must be of type (i), and that
$t\leq \lfloor \frac{n-1}{k-1} \rfloor$ since each step of type (i)
reduces the number of connected components by $k-1$ and the total
reduction of components along the whole chain is $n-1$.

Now, suppose that $1\leq t\leq \lfloor \frac{n-1}{k-1} \rfloor$. A maximal
chain of length  $n-2-t(k-3)$ is constructed as follows. First perform
a sequence of $t$ covering steps of type (i) producing the graph with $k$-cliques
on the sets $\{1,\dots,k\}$, $\{k,\dots,2k-1\}$, $\dots$ 
$\{(t-1)k-(t-2),\dots,tk-(t-1)\}$. Then continue from there via a
sequence of $t-1$ covering steps of type (ii) leading to the graph with
a $(tk-(t-1))$-clique on the set $[tk-(t-1)]$. Finally, $n-(tk-(t-1))$
covering steps of type (iii) will lead to the complete graph. The total
number of steps taken, i.e. the length of the constructed chain, is
$t+(t-1)+n-(tk-(t-1))=n-2-t(k-3)$.
\end{Proof}

The above result yields some nontrivial information about the topology of
$\Delta_{n,k}^2$. For instance, part (i) shows that the order complex of
$\ov{\Sigma_{n,2}}$ is pure of dimension $2n-5$. With Theorem \ref{dnksnk}
this implies that the homology of $\Delta_{n,2}^2$ vanishes in dimensions
greater than $2n-5$ and is free in dimension $2n-5$. Of course, in this
case we already have more precise knowledge from Theorem \ref{main}. By
similar reasoning we can conclude the following new information about
the $k=3$ case from part (ii) of Theorem \ref{chsnk}.

\begin{theorem}
$\widetilde{H}_i(\Delta_{n,3}^2)=0$ for all $i>n-4$, and
$\widetilde{H}_{n-4}(\Delta_{n,3}^2)$ is free.
\label{snk3}
\end{theorem}

In the remaining cases the following can be deduced.

\begin{theorem}
Assume that $k>3$.
\begin{itemize}
\item[(i)] $\widetilde{H}_i(\Delta_{n,k}^2)=0$ if $i>n-k-1$
or $n-k-1>i>n-2k+2$.
\item[(ii)] $\widetilde{H}_{n-k-1}(\Delta_{n,k}^2)$ is free of
dimension ${{n-1} \choose {k-1}}$.
\item[(iii)] If $n<2k-1$ then $\Delta_{n,k}^2$ has the homotopy type
of a wedge of ${{n-1} \choose {k-1}}$ spheres of dimension $n-k-1$.
\item[(iv)] If $n=2k-1$ then $\Delta_{n,k}^2$ has the homotopy type
of a wedge of spheres.  This wedge consists of ${{n-1} \choose {k-1}}$ 
$(n-k-1)$-spheres and $\frac{1}{2}n{{n-1} \choose {k-1}}$ $1$-spheres.
\end{itemize}
\label{snk4}
\end{theorem}

\begin{Proof}
We use Theorem \ref{dnksnk} without reference throughout the proof.
Claim (i) follows immediately from Theorem \ref{chsnk}(iii),(iv),(v).  By
Theorem \ref{chsnk}(v), the subposet of $\snk$ generated by chains of
length $n-k+1$ is isomorphic to the poset obtained by removing all sets
of sizes $1,2,\ldots,k-1$ from the Boolean algebra $B_n$.  Claims (ii)
and (iii) now follow immediately from the rank selection results 
in \cite{Bjo91,Sta82}, along with Theorem \ref{chsnk}(iii),(iv).  

If $n=2k-1$, let $W$ be the set of 
vertices in $\Delta(\ov{\snk})$ corresponding to graphs which consist of two
$k$-cliques intersecting in a single vertex, and let $\Delta_0$ be the
complex obtained by removing all simplices containing an element of $W$
from $\Delta(\ov{\snk})$.  Then $\Delta_0$ is the order complex of the 
subposet of $\ov{\snk}$ generated by chains of length $n-k+1$, and is 
therefore homotopy equivalent to a wedge of ${{n-1} \choose {k-1}}$
$(n-k-1)$-spheres, as above.  If $G \in \snk$ corresponds to an element
$w \in W$, then by Lemma \ref{covsnk}, $(\ov{\snk})_{<G}$ consists of two graphs
which contain a $k$-clique and $n-k$ isolated vertices.  It follows that
$link_{\Delta(\ov{\snk})}(w)$ consists of two vertices in $\Delta_0$.  There
is a homotopy equivalence between $\Delta_0$ and a wedge of
${{n-1} \choose {k-1}}$ $(n-k-1)$-spheres which maps 
$\cup_{w \in W}link_{\Delta(\ov{\snk})}(w)$ to the wedge point.  It is easy
to see that $|W|=\frac{1}{2}n{{n-1} \choose {k-1}}$, and claim (iv) follows.
\end{Proof}

%

\medskip

The homology of $\Da^{2}_{n,3}$ has been computed for $4 \leq n \leq 7$.  
It is concentrated in dimension $n-4$, see Table 2. 

$$
\begin{array}{|c||c|c|c|c|c|c|}
\hline 
n \backslash i & 0 & 1             & 2 & 3 \\
\hline 
\hline 
2 & 0     & 0        &           0 &          0 \\
\hline 
3 &  0    & 0        &           0 &          0 \\
\hline
4 & \ZZ^3 & 0        &           0 &          0 \\
\hline 
5 & 0     & \ZZ^{21} &           0 &          0 \\
\hline
6 & 0     &        0 & \ZZ^{180}   &          0 \\
\hline
7 & 0     &        0 & 0           & \ZZ^{2010} \\
\hline 
\end{array}
$$

\centerline{{\sf Table 2:} Homology groups $\widetilde{H}_i(\Delta_{n,3}^2)$}

We believe that the concentration of homology in dimension $n-4$ is true
in general, see the discussion in Section 9.2. One approach to proving
this could be via the following lemma.
Recall that a graph is called a {\em forest} if it is free of circuits. This is
equivalent to saying that every block in its block decomposition has at
most two vertices.

\begin{lemma}
Suppose that the order complex of the open interval $(G,\hat{1})$ in
$\Sigma_{n,2}$ is topologically $(n-5)$-connected for every forest 
$G$. Then $\Da_{n,3}^2$ is homotopy equivalent to a wedge of
$(n-4)$-spheres.
\label{sn2sn3}
\end{lemma}

\begin{Proof}
By Theorem \ref{dnksnk} we may replace $\Da_{n,3}^2$ by 
$\Delta(\ov{\Sigma_{n,3}})$, which by Theorem \ref{chsnk}(ii) is
$(n-4)$-dimensional. Hence by known reductions (see \cite[(9.19)]{Bjo91})
it suffices to prove that $\Delta(\ov{\Sigma_{n,3}})$ is $(n-5)$-connected.
By Theorem \ref{main} we know that $\ov{\Sigma_{n,2}}$
is $(n-5)$-connected, and we will show how to transfer this connectivity
to the subposet $\ov{\Sigma_{n,3}}$ under the given hypothesis.

Let $P_n$ be the subposet of $\ov{\Sigma_{n,2}}$ consisting of all
elements which contain at least one block of size greater than two.
The elements in $\ov{\Sigma_{n,2}} \setminus P_n$ are the forests $G$,
so a version of Quillen's fiber lemma (see \cite[Lemma 11.12]{Bjo91})
together with our hypothesis about the intervals $(G,\hat{1})$
shows that $P_n$ is $(n-5)$-connected.

Now, note that $\ov{\Sigma_{n,3}} \subseteq P_n$.  
Let $\rho: P_n
\rightarrow \ov{\Sigma_{n,3}}$ be the map which sends $H \in P_n$ to the
subgraph obtained by removing from $H$ all edges which are not contained 
in a block of size at least three. Then $\rho$ is a lower closure 
operator on $P_n$ (that is, a closure operator on $P_n$ with the opposite 
order) whose image is $\ov{\Sigma_{n,3}}$. Hence, by Corollary \ref{closure}
$\ov{\Sigma_{n,3}}$ is $(n-5)$-connected also.
\end{Proof}

We end this section with an easy result which shows that the homology
of $\Da_{n,k}^2$ vanishes in all sufficiently low dimensions.  For this
the posets $\snk$ are not used.

\begin{lemma}
Let $E$ be a $k$-graph on $n$ vertices.  If $E$ is 
$2$-connected, then $E$ contains at least $\lceil\frac{n}{k-1} \rceil$ 
hyperedges.
\label{skellem}
\end{lemma}

\begin{Proof}
If $E$ is $2$-connected then for each $k$-edge $X=\lp v_1,\ldots,v_k \rp
\in E$ there exist at least two $v_i$ which are contained in some $k$-edge
of $E$ other than $X$.  It follows easily that
\[
n \leq |E|(k-1).
\]
\end{Proof}

\begin{corollary}
The complex $\Da_{n,k}^2$ is topologically 
$(\lceil \frac{n}{k-1} \rceil -3)$-connected, implying that
$\widetilde{H}_i(\Delta_{n,k}^2)=0$ for $i=0,1,\dots,
\lceil \frac{n}{k-1} \rceil -3.$
\label{skelcor}
\end{corollary}

\begin{Proof}
Let $m=\lceil \frac{n}{k-1} \rceil -2$.  By Lemma \ref{skellem} 
$\Da_{n,k}^2$ contains the full $m$-skeleton on its vertex set.  The
corollary follows immediately.
\end{Proof}

\section{The Euler characteristic of the complex $\Da_{n,k}^2$} 
\label{notconnhyper}

In Section \ref{closed} we were able to determine the homotopy type of
$\Da_{n,k}^2$ for $k>3$ when $n \leq 2k-1$, but not for $k=3$, nor for 
$k>3$ and $n>2k-1$.  Indeed, other than the connectivity result given
in Corollary \ref{skelcor}, in the case $k=3$ our only information on
the topology of $\Da_{n,k}^2$ is given by Theorem \ref{snk3} 
(unless $n$ is very small), and in the case
$k>3$ and $n>2k-1$ we were only able to determine
the homology group $\widetilde{H}_{n-k-1}(\Da_{n,k}^2)$.  

In this section,
we investigate the reduced Euler characteristic of $\Da_{n,k}^2$. 
We will determine a formula for the exponential generating 
function
$$M_k(x):=\sum_{n=1}^\infty \tilde{\chi}(\Da^{2}_{n,k})\frac{x^n}{n!},$$
for all $k \geq 2$.  That formula is stated in the following theorem.

\begin{theorem}\label{th:mobius}
For $k\ge 2$, we have
\[M_k'\left(x\frac{p_{k-1}(x)}{p_k(x)}\right)=
\ln\left(\frac{p_{k-1}(x)}{p_k(x)}\right),
\]
where $p_k(x):=1+x+\frac{x^2}{2!}+\cdots +\frac{x^{k-1}}{(k-1)!}$.
\end{theorem}

Theorem \ref{th:mobius} gives another proof that 
$\tilde{\chi}(\Da^{2}_{n,2})=-(n-2)!$.  It also implies 
\[M_3'(x)=\ln\left(\frac{-x(x-2)}{(x-1)+\sqrt{2-(x-1)^2}}\right),
\]
which gives the sequence $0,0,-1,3,-21,180,-2010,27090,-430290,\ldots$
for $\tilde \chi(\Da^2_{n,3})$, cf. Table 2.
To obtain these corollaries set $y:=x\frac{p_{k-1}(x)}{p_k(x)}$ and solve
for $x$ to get $x=\frac{1}{1-y}$ and $x=\frac{y-1+\sqrt{2-(y-1)^2}}{2-y},$
when $k=2$ and $3$ respectively. 

\medskip

To prove Theorem \ref{th:mobius} we will use the posets $\snk$ defined
in Section \ref{closed}.  Note that 
$\tilde\chi(\Da^2_{n,k})=\mu_{\Sigma_{n,k}}(\hat 0,\hat 1)$.
We will write $\mu_k(G)$ for $\mu_{\Sigma_{n,k}}(\hat 0,G)$ and $\mu_k(n)$ 
for $\mu_{\Sigma_{n,k}}(\hat 0,\hat 1)$.

Let $\Pi_{n,k}$ be the $k$-equal lattice, which is the lattice of
partitions of $\lb n \rb$ into subsets such that each subset
has size one or at least $k$.  Let 
$\tau_k(n):=\mu_{\Pi_{n,k}}(\hat 0,\hat 1)$ for $n\ge k$ and 
$\tau_k(2)=\cdots=\tau_k(k-1)=0$, but $\tau_k(1)=1$.  
The exponential generating 
function, $T_k(x):=\sum_{n=1}^\infty \tau_k(n)\frac{x^n}{n!}$, 
for the M\"obius function of $\Pi_{n,k}$ is known to be 
\[ T_k(x)=\ln(p_k(x)), 
\] 
where $p_k(x)$ is as above.  It was first
calculated in \cite{BL}.  

Let $\cC$ be the set of connected graphs in $\Sigma_{n,k}$.
Now let
\[ \sigma: \Sigma_{n,k}\setminus \cC \longrightarrow \Pi_{n,k}
\]
be the function which maps a disconnected graph in $\snk$ to the partition 
determined by 
its connected components.  It is easily seen that for each $x \in \Pi_{n,k}$,
$\sigma^{-1}_{\leq}(x)$ has a unique maximum element and therefore has a
contractible order complex.  Thus by Proposition \ref{quillen} and the
definition of the M\"obius function, we have  
\[\tau_k(n)=-\sum_{G\in \Sigma_{n,k}\setminus \cC} \mu_k(G)=
\sum_{G\in \cC} \mu_k(G).
\]
Thus it suffices to concentrate on the connected but not $2$-connected graphs.
First we need a simple lemma.
\begin{lemma}
If $G \in \cC$ has blocks $W_1,\ldots,W_r$ with $|W_i|=w_i$, then
\[\mu_k(G)=\prod_{i=1}^{r}\mu_k(w_i).
\]
\label{mult}
\end{lemma}
\begin{Proof} By Lemma \ref{interv} the 
interval $\lb \hat{0},G \rb$ is isomorphic to the product poset
$\Sigma_{w_1,k} \times \ldots \times \Sigma_{w_r,k}$.  The lemma now
follows from the well known multiplicativity of the M\"obius function.
\end{Proof}

Now define
$\displaystyle{\al_k(n):=\sum_{G\in\cC_1} \mu_k(G)}$, where
$\cC_1:=\{G\in \cC: \mbox{$n$ is not a cutpoint of G}\}$.  Also, set
$\al_k(1)=\cdots=\al_k(k-1)=0$ and 
$A_k(x):=\displaystyle{\sum_{n=1}^\infty \al_k(n)\frac{x^n}{n!}}$. 

\begin{lemma}\label{lm:Dk}
We have
\[
A_k'(x)=\ln\left(\frac{p_{k-1}(x)}{p_k(x)}\right).
\]
\end{lemma}

\begin{Proof}
If $n$ is a cutpoint of $G \in \cC$, let $P_1,\dots,P_t$ be the connected
components of $G-n$.  Then for each $i \in \lb t \rb$, $n$ is not a 
cutpoint of the connected subgraph of $G$ induced on $P_i \cup \lp n \rp$.
Using Lemma \ref{mult}, we get the recursive formula
\[\tau_k(n)=\al_k(n)+\sum_{t=2}^{\left\lfloor\frac{n-1}{k-1}\right\rfloor}
\sum_{P_1|\ldots|P_t \in \Pi_{n-1}}\al_k(|P_1|+1)\cdots\al_k(|P_t|+1),
\]
where each summand in the double sum on the right counts the M\"obius 
functions of all elements $G \in \cC$ such that $G-n$ has connected 
components $P_1,\ldots,P_t$.  By the definition of $\al_k(n)$ we can
rewrite this formula as
\[
\tau_k(n)=\sum_{P_1|\ldots|P_t \in \Pi_{n-1}}
\al_k(|P_1|+1)\cdots\al_k(|P_t|+1).
\]
The exponential formula (Proposition \ref{exp}) and easy 
power series manipulations then give
\[T_k'(x)=e^{A_k'(x)}.
\]
\end{Proof}

\noindent
{\sf Proof of Theorem \ref{th:mobius}:}
We will establish a recurrence relation for $\al_k$ 
involving $\mu_k$.  Let $G \in \cC_1$ and let $W$ be the block of $G$
containing $n$.  Then $W$ is the unique maximal clique in $G$ which 
contains $n$.  Let $S=W \setminus \lp n \rp$.  Let $B_G$ be the graph on
the blocks of $G$ defined in Proposition \ref{sep}.  Let $T_1,\ldots,T_t$
be the connected components of $B_G-W$ and for each $T_i$ let $V_i$ be the
set of vertices of $G$ which are contained in a block that is contained
in $T_i$.  For each $T_i$ there is a unique $j_i \in S$ such that the
subgraph of $G$ induced on $V_i \cup \lp j_i \rp$ is a connected union of 
blocks of $G$.  Conversely, any $G \in \cC_1$ can be obtained by choosing $S$, 
$T_i$ and $j_i$ as above and then choosing graphs $H_i$ on $T_i \cup 
\lp j_i \rp$ such that each $H_i$ is either a clique of size at least $k$ or 
isomorphic to a connected element of $\ov{\Sigma}_{|T_i|+1,k}$.  These 
choices can all be made independently, so we get the recurrence relation
\[
\al_k(n)=\sum_{S\uplus T=[n-1]}\mu_k(|S|+1)\sum_{T_1\uplus\ldots\uplus T_t=T}
(\al_k(|T_1|+1)|S|)\cdots(\al_k(|T_t|+1)|S|).
\]
We cannot apply the exponential formula directly at this point due to the factors
$|S|$ which appear on the right hand side of the above equation.  However, we
get    

\[A_k'(x)=\sum_{n=1}^\infty \frac{\al_k(n)x^{n-1}}{(n-1)!}
\]
\[=\sum_{n=1}^\infty\sum_{i=1}^{n-1}\binom{n-1}{i}\frac{\mu_k(i+1)x^i}{(n-1)!}
\sum_{T_1|\ldots|T_t \in \Pi_{n-i-1}}(\al_k(|T_1|+1)i)\cdots(\al_k(|T_t|+1)i)
x^{n-i-1}
\]
\[=\sum_{i=1}^\infty\frac{\mu_k(i+1)x^i}{i!}\sum_{n=i+1}^\infty
\sum_{T_1|\ldots|T_t \in \Pi_{n-i-1}}
(\al_k(|T_1|+1)i)\cdots(\al_k(|T_t|+1)i)\frac{x^{n-i-1}}{(n-i-1)!}
\]
Applying the exponential formula, for each $i$ we get
\[
\sum_{n=i+1}^\infty
\sum_{T_1|\ldots|T_t \in \Pi_{n-i-1}}
(\al_k(|T_1|+1)i)\cdots(\al_k(|T_t|+1)i)\frac{x^{n-i-1}}{(n-i-1)!}
=e^{iA_k'(x)}.
\]
Thus
\[A_k'(x)=\sum_{i=1}^\infty\frac{\mu_k(i+1)x^i}{i!}e^{iA_k'(x)}=
M_k'\left(xe^{A_k'(x)}\right).
\]
The theorem now follows from Lemma \ref{lm:Dk}.
\qed

%
%
%

\section{$(n-2)$-connected graphs and matching complexes} \label{match}

Before we proceed to consider $(n-2)$-connected graphs, let us 
state some simple but useful facts about the general situation. 
What do maximal $(i-1)$-separable
graphs on the $n$ element set $[n]$ look like?  Is is clear
that each such graph is described by an $(i-1)$-set $A$ and a partition 
$B \uplus C$ of $[n] \setminus A$ into two non-empty 
blocks $B$, $C$. 
The corresponding maximal $(i-1)$-separable 
graph is the complete graph on $[n]$ with all edges 
connecting $B$ and $C$ removed.

Now let $G$ be an $(n-2)$-connected graph on $n$ vertices, so $G \not\in
\Delta_{n}^{n-2}$. Then by 
the above description of maximal $(n-3)$-separable graphs the induced
subgraph on any three vertices must contain at least two edges. Thus
the complementary graph (i.e., the graph containing 
precisely the edges that are
not in $G$) is a matching. The graphs on $n$ vertices that are matchings form a
simplicial complex, that we denote by $M_n$. We conclude the following.

\begin{proposition} The matching complex $M_n$ is Alexander dual (in the sense
of Proposition \ref{alexander}) to the complex $\Delta_n^{n-2}$. In
particular, there is an isomorphism  
$$\widetilde{H}_i(M_n) \cong \widetilde{H}^{\binom{n}{2}-i-3}(\Delta_n^{n-2}).$$
\end{proposition}

The matching complexes $M_n$ have attracted attention for
various reasons. 
In \cite[Theorem 4.1]{BLVZ92} the matching complex $M_n$ is
shown to be topologically $( \lfloor \frac{n+1}{3} \rfloor-2)$-connected, 
which implies that $\widetilde {H}_i(M_n)=0$,
for $i=0,\dots, \lfloor \frac{n+1}{3} \rfloor -2$.  
We thus get the following corollary.

\begin{corollary} 
The cohomology of $\Delta_n^{n-2}$ vanishes in
dimensions $i \geq \binom{n}{2}-\lfloor \frac{n+1}{3} \rfloor -1$.
\end{corollary}

The following table shows what we know about the homology groups 
$\widetilde{H}_i(M_n)$, based on the results of \cite{BLVZ92} for $n\leq 6$ 
and $n=8$, and our own
computations.  

$$
\begin{array}{|c||c|c|c|c|c|c|}
\hline 
n \backslash i & 0 & 1       & 2 & 3 & 4 & 5 \\
\hline 
\hline 
2 & 0     & 0        & 0 & 0 & 0 & 0 \\
\hline 
3 &  \ZZ^2 & 0        & 0 & 0 & 0 & 0 \\
\hline
4 & \ZZ^2 & 0        & 0 & 0 & 0 & 0 \\
\hline 
5 & 0     & \ZZ^6                     & 0 & 0 & 0 & 0 \\
\hline
6 & 0     & \ZZ^{16}                  & 0 & 0 & 0 & 0 \\
\hline
7 & 0     & \text{torsion}^1& \ZZ^{20}  & 0 & 0 & 0 \\
\hline 
8 & 0     & 0        & \ZZ^{132} & 0 & 0 & 0 \\
\hline 
9 & 0     & 0        & \ZZ^{42} \oplus \text{torsion}^2 
& \ZZ^{70} & 0 & 0 \\
\hline
10 & 0     & 0        & \text{torsion}^3 & \ZZ^{1216} & 0 & 0\\
\hline
11 & 0     & 0        &                0 & \ZZ^{1188} \oplus \text{torsion}^4 & \ZZ^{252} & 0\\
\hline
12 & 0     & 0        &                0 & \text{torsion}^5 & \ZZ^{12440} & 0\\
\hline
\end{array}
$$

\centerline{{\sf Table 3:} Homology groups $\widetilde{H}_i(M_n)$ of matching complexes}

\medskip
\footnotetext[1]{There is $\ZZ_3$-torsion of rank $1$. No $\ZZ_p$-torsion 
for $p=2$, $5 \leq p \leq 17$.}
\footnotetext[2]{There is $\ZZ_3$-torsion of rank $8$. No $\ZZ_p$-torsion
for $p=2$, $5 \leq p \leq 17$.}
\footnotetext[3]{There is $\ZZ_3$-torsion of rank $1$. No $\ZZ_p$-torsion
for $p=2, 5, 7$.}
\footnotetext[4]{There is $\ZZ_3$-torsion of rank $35$. No $\ZZ_p$-torsion
for $p=2, 5, 7$.}
\footnotetext[5]{There is $\ZZ_3$-torsion of rank $56$. No $\ZZ_p$-torsion
for $p=2, 5, 7$.}
We see that the complexes  $\Delta_n^{n-2}$ can have torsion, 
and that this phenomenon begins with $\Delta_{7}^{5}$.

\section{The Euler characteristic of the complex $\Da_{n}^{n-3}$} 
\label{notn3}

Consider the complex $(\Delta_n^{n-3})^*$ which is 
the Alexander dual of $\Delta_n^{n-3}$,
and the exponential generating function of
its reduced Euler characteristic $$F_n^{n-3} (x) := 
\displaystyle{\sum_{n \geq 0} \widetilde{\chi}((\Delta_n^{n-3})^*) \frac{x^n}{n!}}.$$ 
The values of $\widetilde{\chi}((\Delta_n^{n-3})^*)$ in the degenerate cases
$n\leq 3$ will appear from an explicitly calculated expansion below.
We will express the reduced Euler characteristic of $\Delta_n^{n-3}$ in
terms of an expression for this series.

\begin{theorem}  
We have that:
$$F_n^{n-3}(x)=
x-\frac{\exp(\frac{x}{2(1+x)})+x-\frac{1}{4}x^2-\frac{1}{8}x^4}{\sqrt{1+x}}.$$
The exponential generating function of the reduced Euler characteristic of
$\Delta_n^{n-3}$ is then the sum of the real and imaginary parts of
$-F_n^{n-3}(ix)$.
\end{theorem}
\Proof We will argue as we did for $\Delta_n^{n-2}$ in Section \ref{match}.
If a graph $G$ is $(n-3)$-connected then the induced subgraph on any $4$ of
its vertices contains either a vertex of degree $3$ or a path of length $3$. 
Thus in the complementary graph the induced subgraph on any $4$ vertices is
either contained in a $3$-cycle or in a path of length $3$. In particular, 
there are no $4$-cycles and no vertices of degree $3$ in the complementary 
graph. Thus the connected components of the complementary graph  are paths of 
any length and cycles of length different from $4$. Moreover, any graph in 
which every connected component is a cycle not of length $4$ or a path is the
complement of an $(n-3)$-connected graph, so $(\Delta_{n}^{n-3})^*$ consists
of all such graphs.  There are exactly $n!/2$ different paths of 
length $n$ and $(n-1)!$ different $n$-cycles on an $n$-element vertex
set. Now a direct application of the Exponential Formula \ref{exp}
gives the result for the generating function of $(\Delta_n^{n-3})^*$. The 
remaining assertion follows from the fact that when passing from $\Delta_n^i$ 
to its Alexander dual the Euler characteristic changes by a factor of 
$-(-1)^{n(n-1)/2}$.
\qed

The complex $(\Delta_{n}^{n-3})^*$ has maximal simplexes of dimensions 
$n-1$ and $n-2$ only.
It is easily collapsible to a pure complex of dimension $n-2$.
A Maple computation (see below) shows that neither the Euler characteristic
of $(\Delta_n^{n-3})^*$ nor the Euler characteristic of $\Delta_n^{n-3}$
alternate in sign, so the pure complexes are certainly not all 
Cohen-Macaulay.  The calculation shows that

                      
$$F_n^{n-3}(x)=-1-x+\frac{1}{4}x^4+\frac{1}{20}x^5+\frac{1}{20}x^6
-\frac{1}{27}x^7-
\frac{1}{224}x^8-\frac{1}{480}x^9+O(x^{10}).$$

%
%
%
%

We have also studied the slightly larger complex of graphs which are
the disjoint union of cycles and paths of any lengths (i.e., graphs 
with maximum vertex degree at most $2$). This is also a reasonable
generalization of the matching complex, which is the complex of all 
graphs with maximum vertex degree at most $1$.   
The Euler characteristic of the corresponding Alexander dual 
has almost the same generating function
as for $(\Delta_n^{n-3})^*$.  That generating function is 

$$x - \frac{\exp(\frac{x}{2(1+x)})+x-\frac{1}{4}x^2}{\sqrt{1+x}} = $$

$$-1-x+\frac{1}{8}x^4-\frac{3}{40}x^5+\frac{1}{20}x^6-\frac{1}{28}x^7-
\frac{17}{896}x^8-\frac{7}{1920}x^9-\frac{23}{2400}x^{10}+ O(x^{11}).$$

The maximal simplices in the cycles-and-paths complex have dimension $n-1$
or $n-2$, and the complex can be collapsed to a pure $(n-2)$-dimensional
complex.  The generating function for the Euler characteristic shows
that these collapsed complexes are not all Cohen-Macaulay.

%
%
%
%
%
%
%
%

\section{Final Remarks}

\subsection{Homology and Topology of $\Delta_{n}^i$} 

The results and computations presented in this paper suggest that there 
is probably no uniform statement that covers the topology of all complexes
$\Delta_n^i$. However, for $i \leq 2$ the homotopy type
calculations for $\Delta_n^i$  give very nice answers.
This is consistent with the graph theoretical study
of not $i$-connected graphs, where there is a
good structure theory only when $i \leq 3$ (see for example Chapter 6 of 
Lov\'asz' book \cite{Lov93}, or the survey article by Oxley \cite{Oxl96} and
the references therein). 

As mentioned, there is a structure theory for $3$-connected graphs. The $3$-connected 
graphs on $n$ vertices for which neither the deletion nor the contraction of
an edge leads to a $3$-connected graph were classified by Tutte 
(see Theorem 2.3 in \cite{Oxl96}) as ``wheels'' and ``whirls,''
both having $2n-2$ edges. Note that
this does not provide a characterization of the deletion-minimally $3$-connected 
graphs (the minimal non-faces of $\Delta_n^3$),
however it does show that no graph with less than $2n-2$ edges can be
$3$-connected. Hence, $\Delta_n^3$ has a complete $2n-4$-skeleton,
which shows that $\widetilde{H}_i(\Delta_n^3)=0$ for $i<2n-4$.  This
fact together with the following table leads to an interesting conjecture.

$$
\begin{array}{|c||c|c|c|c|c|c|c|c|c|c|c|}
\hline 
n \backslash i & 0 & 1 & 2 & 3 &     4 & 5 &     6 & 7 &        8 & 9 & 10 \\
\hline 
\hline 
3              & 0 & 0 & 0 & 0 &     0 & 0 &     0 & 0 & 0        & 0 & 0\\
\hline
4              & 0 & 0 & 0 & 0 & \ZZ^1 & 0 &     0 & 0 & 0        & 0 & 0\\
\hline 
5              & 0 & 0 & 0 & 0 &     0 & 0 & \ZZ^6 & 0 & 0        & 0 & 0\\
\hline
6              & 0 & 0 & 0 & 0 &     0 & 0 &     0 & 0 & \ZZ^{36} & 0 & 0\\
\hline
7              & 0 & 0 & 0 & 0 &     0 & 0 &     0 & 0 &        0 & 0 & \ZZ^{240} \\
\hline 
\end{array}
$$

\centerline{{\sf Table 4:} Homology groups $\widetilde{H}_i(\Delta_n^3)$}

\medskip

\noindent {\sf Conjecture 9.1:} $\Da_n^3$ has the homotopy type of a wedge
of $\frac{(n-3)(n-2)!}{2}$ spheres of dimension $2n-4$.

\medskip

For general $i$ the situation (concerning $\Da_n^i$)
seems to be far more complicated. 
However, we would like to remark that for $i=2,3$ the known Betti
numbers are (up to sign) the Lah-numbers $L_{n-2,1}$ and $L_{n-2,2}$
(see \cite[p.165-166]{Com70}). This coincidence unfortunately fails for $i=4$,
which is easily seen by comparing $L_{n-2,3}$ with $\widetilde{\chi}(\Delta_n^4)$ for $n=6$. For $i>3$ no good structure theory for
$i$-connected graphs is known, and the results 
of Section \ref{match} on $(n-2)$-connected
graphs indicate that the topology of $\Delta_n^i$ will not behave nicely
for all $i$.
Nevertheless, the Alexander Duality with the matching complexes $M_n$
encourages a closer look at the complexes $\Delta_n^{n-2}$. Surprisingly, 
the prime $3$ seems to play a special role in the topology of these
complexes. We have not detected $p$-torsion
in the homology of $\Delta_n^{n-2}$ for any prime $p \neq 3$.

\medskip

\noindent {\sf Question 9.2:} Does the homology of $M_n$ have $p$-torsion
for any prime $p \neq 3$?

\medskip

Even more surprising is the fact that the same prime $3$ seems to
play an analogous role for the matching complexes $M_{n,n}$ on complete
bipartite graphs $K_{n,n}$, also called chessboard complexes
(see \cite{BLVZ92}).

$$
\begin{array}{|c||c|c|c|c|c|c|c|}
\hline 
n \backslash i & 0 & 1 & 2 & 3 & 4 & 5 & 6 \\
\hline 
\hline 
1              & 0   & 0     & 0          &            0 &        0 & 0 & 0 \\
\hline
2              & \ZZ & 0     & 0          &            0 &        0  & 0 & 0\\
\hline
3              & 0   & \ZZ^4 & 0          &            0 &        0 & 0 & 0 \\
\hline
4              & 0   &     0 & \ZZ^{15}   &            0 &        0 & 0 & 0 \\
\hline
5              & 0   &     0 &  \ZZ_3     & \ZZ^{56}   &          0 & 0 & 0 \\
\hline 
6              & 0   &     0 &          0  & \ZZ^{25}\oplus \text{torsion}^1 & \ZZ^{210}   & 0  &0 \\
\hline
7              & 0   &     0 &          0 &          0 & \ZZ^{588} \oplus \text{torsion}^2 & \ZZ^{792} & 0 \\
\hline 
8              & 0   &     0 &          0 &          0 & \text{torsion}^3 & ? & ? \\
\hline 
\end{array}
$$

\footnotetext[1]{There is $\ZZ_3$-torsion of rank $10$. No $\ZZ_p$-torsion for $p=2,5,7$}
\footnotetext[2]{There is $\ZZ_3$-torsion of rank $66$. No $\ZZ_p$-torsion for $p=2,5,7$}
\footnotetext[3]{There is $\ZZ_3$-torsion of rank $1$. This group is finite
according to \cite{FriHan96}.}

\centerline{{\sf Table 5:} Homology groups $\widetilde{H}_i(M_{n,n})$ for the bipartite matching complexes}

\medskip
It is easy to see that $M_{n,n}$ collapses to an $(n-2)$-dimensional
complex. Hence homology is free in dimension $n-2$ and vanishes in higher
dimensions. 
Looking at the table and the footnotes the following question naturally occurs.

\medskip

\noindent {\sf Question 9.3:} Is the homology of $M_{n,n}$ 
free except for $3$-torsion ?

\medskip

\subsection{Homology and Topology of $\Delta_{n,k}^2$} 
For $k$-graph complexes the problem of determining the topology of the
complex of separable $k$-graphs (the $i=2$ case) seems to be the 
most important. The complexes $\Delta_{n,k}^2$ play the 
same role in the study 
of spaces of ``knots'' for which $k$-fold self-intersections are forbidden as 
the complexes $\Delta_n^2$ play for ordinary knots.

\medskip

\noindent {\sf Question 9.4:} What is the homology and homotopy type of
$\Delta_{n,k}^2$ ? 

\medskip

The evidence from Section \ref{closed} leads us to anticipate the
following answer for $k=3$.

\medskip

\noindent {\sf Conjecture 9.5:} 
$\Da_{n,3}^2$ is homotopy equivalent to a wedge of
$(n-4)$-spheres.
\medskip

A natural approach to this question is through further
combinatorial study of the lattices
$\snk$ defined in Section \ref{closed}.

\medskip

\noindent {\sf Conjecture 9.6:} The lattice $\Sigma_{n,k}$ is shellable.

\medskip

This is open in all cases, except for the somewhat degenerate cases $n<2k-1$
when $\snk$ is a truncated Boolean algebra. If Conjecture 9.6 were
verified for $k=2$ it would reprove Theorem \ref{main}, and it would
via Lemma \ref{sn2sn3} imply the truth of Conjecture 9.5. If Conjecture
9.6 were verified for $k=3$ it would also imply the truth of Conjecture 9.5.
If Conjecture 9.6 were verified for $k>3$ it would via Theorem \ref{chsnk}(iv)
and the results of \cite{BjoWac96} imply the truth of the following.

\medskip

\noindent {\sf Conjecture 9.7:} If $k>3$ and $n \geq 2k-1$, then
$\Da_{n,k}^2$ is homotopy equivalent to a wedge of spheres. Furthermore,
the dimensions of the spheres are precisely $n-4-t(k-3)$
for $1\leq t\leq \lfloor \frac{n-1}{k-1} \rfloor$.

\medskip

\subsection{Generating series of Euler characteristics} 

Let $F_i(x) = \displaystyle{\sum_{n \geq 0} \widetilde{\chi}(\Delta_{n}^{i})
\frac{x^n}{n!}}$
and $G_i(x) = \displaystyle{\sum_{n \geq 0} \widetilde{\chi}((\Delta_{n}^{n-i})^*)
\frac{x^n}{n!}}$.
By the results presented in this paper we get the following table:

\medskip

\begin{center}
$\begin{array}{|c||c|c|}
\hline
i & F_i(x) & G_i(x) \\
\hline 
\hline
1 & ln(1+x) & -1 \\
\hline
2 & (1-x)\log(1-x)+1+x &  -\exp(x-\frac{x^2}{2}) \\
\hline 
3 &                   & x - \frac{\exp(\frac{x}{2(1+x)})+x-\frac{1}{4}x^2-\frac{1}{8}x^4}{\sqrt{1+x}} \\
\hline
\end{array}$
\end{center}

\centerline{{\sf Table 6:} Generating functions of the Euler characteristics of $\Delta_n^i$ and $(\Delta_n^{n-i})^*$}

\medskip

We cannot formulate a conjecture about the entries in this table
for $i>3$. Nevertheless, even though the
actual homology computation may be too difficult, the 
generating series may be computable for a few more cases.
Assuming a positive answer to Conjecture 9.1, we get

$$F_3(x) = (x-\frac{3}{2})\log(1-x)+1-\frac{3}{2}x+\frac{1}{4}x^2.$$

\subsection{The representation of the symmetric group} \label{questionrep} 

All complexes $\Delta_{n,k}^i$ are invariant under the action of the
symmetric group $S_n$. This action determines a linear representation
of $S_n$ on each homology group $\widetilde{H}_j(\Delta_{n,k}^i)$. 
For fixed $n,k$ and $i$, the alternating sum of the characters of the
given representations is a virtual character of $S_n$ 
that we denote by $\omega_{n,k}^i$.
For $k=2$ and $i=1,2$ 
this is an actual character (up to sign) and satisfies 
$$(-1)^{n+1}\omega_{n+1,k}^2  = (sign_n \omega_{n,k}^1) \uparrow_{S_n}^{S_{n+1}} 
+ sign_{n+1} \omega_{n+1,k}^1.$$

{From} looking at the dimensions of the homology modules it is clear that
the analogous formula for $k \geq 3$ does not hold. We have also seen that
homology of $\Delta_{n,k}^i$ is not torsion-free in general. On the
other hand, \cite{ReiRob97} has demonstrated that for the closely related
matching complexes it is possible to determine the representations
on the rational homology. Thus it is reasonable to ask:

\medskip

\noindent {\sf Question 9.8:} What is the character of $S_n$ on each 
non-vanishing rational homology group of $\Delta_{n,k}^i$ ?

\medskip

The character $\omega_n^2 = \omega_{n,2}^2$ determined in Section \ref{homrep}
has recently appeared in several different areas of mathematics.
First in the work of  C. A. Robinson \& S. Whitehouse \cite{RobWhi94} and
S. Whitehouse \cite{Whi94-2} on gamma-homology of algebras and later in work
of E. Getzler \& M. Kapranov \cite{GetKap95} on operads,
O. Mathieu \cite{Mat96} on hyperplane arrangements and symplectic geometry, 
in the work of M. Kontsevich \cite{Kon93} on Lie algebras and symplectic 
geometry, and in the work of P. Hanlon \cite{Han96}, P. Hanlon \& R.P. 
Stanley \cite{HanSta95} and S. Sundaram \cite{Sun96} in a combinatorial 
and representation-theoretic context. It seems mysterious that the same 
character pops up in so many seemingly unrelated places.

\medskip

\noindent {\sf Question 9.9:} What are the deeper connections between the
various contexts where the character $\omega_n^2$ appears ?

\medskip

The analogous question for the character $\omega_{n,2}^1=
sign_n \cdot \lie_n$ has been
studied quite extensively (see for example \cite{Bar90,BarBer90,Reu93}),
and for that case much detailed information is known. An important aspect
of the work in \cite{Bar90} and \cite{BarBer90} is the 
construction of explicit bases for the modules under consideration.
Thus a first step towards an answer to Question 9.9 could be a solution
of the following problem.

\medskip

\noindent {\sf Problem 9.10:} Describe a combinatorial basis for the
homology of $\Delta_n^2$.

\medskip

A positive answer to the shellability conjecture 9.6 for $k=2$ could
via the induced shelling basis (see \cite{BjoWac96})
lead to progress on Problem 9.10.
 
\section{Notation and Tools}\label{notation}

In this short section we will summarize the main tools that we use in 
the study of the complexes $\Delta_{n,k}^i$. We refer the reader to 
the survey paper \cite{Bjo91} for more details and references.

Let $P$ be a finite partially ordered set -- {\it poset} for short. 
If $P$ has a unique
minimum element $\hat{0}$ and a unique maximum element $\hat{1}$, we denote
by $\ov{P}$ the {\it proper part} of $P$, that is the  
poset obtained by removing from $P$ the elements $\hat{0}$ and $\hat{1}$.
By $\Delta(P)$ we denote the simplicial complex of all chains in ${P}$.
The complex $\Delta(P)$ is called the {\it order complex} of $P$.

By convention we include the empty set $\emptyset$ in every simplicial complex.
For any simplicial complex $\de$, the {\it face lattice} 
$\Lat(\de)$ is the poset of faces
of $\de$, ordered by inclusion and enlarged by an additional greatest element
$\hat{1}$. Then the order complex $\Da(\ov{\Lat(\Da)})$ of the proper 
part of $\Lat(\Delta)$ is homeomorphic to $\Delta$. Indeed,
$\Delta(\ov{\Lat(\Delta)})$ is the barycentric subdivision of $\Delta$.

For a poset $P$ and $p \in P$ we denote by $P_{\leq p}$ the sub-poset 
$\{ p'~|~p' \in P;~p' \leq p~\}$. The posets
$P_{\geq p}$, $P_{< p}$ and $P_{> p}$ are analogously defined. For $p \leq p'$ 
in $P$ we denote
by $[p,p']$ the {\it closed interval} $P_{\geq  p} \cap P_{\leq p'}$ in $P$, and
by $(p,p')$ the {\it open interval} $P_{>p} \cap P_{<p'}$.

For a poset $P$ we denote by $\mu_P$ the $\ZZ$-valued {\it M\"obius function} 
(see \cite{Sta86}), defined recursively on the intervals of $P$ by $\mu_P(x,x) = 1$ and
$\mu_P(x,y) = \displaystyle{-\sum_{x \leq z < y} \mu_P(x,z)}$ if
$x < y$.

By a map $f : P \rightarrow Q$ of posets we always mean a
poset homomorphism (i.e., $x \leq y$ implies $f(x) \leq f(y)$).
For an element $q \in Q$ we denote by $f^{-1}_{\leq}(q)$ the preimage
of $Q_{\leq q}$ under $f$. The poset $f^{-1}_{\geq}(q)$ is analogously defined.

\begin{proposition}[Quillen Fiber Lemma {\rm \cite{Qui78}}] \label{fibre}
Let $f : P \rightarrow Q$ be a map of posets. If $\Delta(f^{-1}_{\leq}(q))$ is
contractible for all $q \in Q$ then $\Delta(P)$ and $\Delta(Q)$ are homotopy
equivalent. 
\label{quillen}
\end{proposition}

A map $f : P \rightarrow P$ from a poset to itself is called a {\it closure
operator} if $f(x) \geq x$ and $f(f(x)) = f(x)$ for all $x\in P$.
The Quillen Fiber Lemma immediately implies the fact that closure operators
preserve the homotopy type.

\begin{corollary}[Closure Lemma]
 \label{closure}Let $f : P \rightarrow P$ be a closure operator on the partially ordered set $P$.
Then $\Delta(P)$ and $\Delta(f(P))$ are homotopy equivalent.
\end{corollary}

If the poset $P$ is a {\it lattice} (i.e., suprema, denoted by
``$\vee$'', and infima, denoted by ``$\wedge$'',  exist) 
then there is another tool for computing the homotopy type. 
Note that if $P$ is a finite lattice then there is a least element 
$\hat{0}$ and a largest element $\hat{1}$ in $P$. For an
arbitrary element $p \in P$ we say that $a \in P$ is a {\it complement} of $p$
if $p \wedge a = \hat{0}$ and $p \vee a = \hat{1}$.

\begin{proposition}[Homotopy Complementation Formula {\rm \cite{BjoWal83}}~] 
\label{homcom}
{\ }
\begin{itemize}
\item[(i)] Let $P$ be a poset and $A \subseteq P$ an
antichain. Assume $\Delta(P \setminus A)$ is contractible. Then 
$\Delta(P)$ is homotopy equivalent to 
$$\bigvee_{x \in A} \Sigma \Bigl(\Delta(P_{<x}) * \Delta(P_{>x}) \Bigr).$$
\item[(ii)] Let $P$ be the proper part of
a lattice and let $Co$ be the set of complements of some
element $p \neq \hat{0}, \hat{1}$. Then $\Delta(P \setminus Co)$
is contractible.
\end{itemize}
\end{proposition}

In the formulation of the proposition $\bigvee$ denotes the wedge
product, $\Sigma$ denotes the suspension and $*$ denotes the join of 
topological spaces.
 
 Our next tool is the combinatorial version of a standard duality theorem from
 algebraic topology.

\begin{proposition}[Combinatorial Alexander Duality]
 \label{alexander} Let $\Delta$ be a finite simplicial complex on vertex set 
$V$ and define
$$\Delta^* = \{ B \subseteq V~|~V \setminus B \not\in \Delta \}.$$
Then
$$\widetilde{H}_i(\Delta) \cong \widetilde{H}^{|V|-i-3}(\Delta^*).$$
\end{proposition}

This is derived as follows. The usual Alexander duality theorem (see e.g.
Munkres \cite{Mun84}) says that $$\widetilde{H}_i(A) \cong \widetilde{H}^{n-i-1}
(S^n \setminus A)$$ for any compact subset $A$ of the $n$-sphere
$S^n$. In our situation, let $P = 2^V \setminus \{ \emptyset, V \}$. This
truncated Boolean algebra is the proper part of the 
face lattice of the boundary complex of
a simplex, so $\Delta(P) \cong S^{|V|-2}$. Now let $A$ be the realization 
of $\Da(\ov{\Lat(\Da)})$ as a subspace of $\Delta(P)$. It is easy to see that
$\Da(P\setminus \ov{\Lat(\Da)})$ is a strong deformation retract
of $S^{|V|-2} \setminus A$, and since $P\setminus \ov{\Lat(\Da)} \cong
\ov{\Lat(\Da^*)}$ the result follows.

Finally, we recall a result from enumerative combinatorics. 
For a number sequence $(a_n)_{n \geq 0}$
the formal power series
$\displaystyle{\sum_{n\geq 0} a_n \frac{x^n}{n!}}$ is called its
{\it exponential generating function}.

\begin{proposition}[Exponential formula] \label{exp}
Suppose that two functions $a,b:\NN \longrightarrow \ZZ$ are given such that 
\[b(n)=\sum_{S_1|\cdots| S_t \in \Pi_n}a(|S_1|)\cdots a(|S_t|),\quad n\ge 1,
\]
where the sum ranges over all set partitions
of $[n]$ and $a(0)=0$, $b(0)=1$. 
Then the exponential generating functions 
$A(x):=\sum_{n=0}^\infty \frac{a(n)x^n}{n!}$ and 
$B(x):=\sum_{n=0}^\infty \frac{b(n)x^n}{n!}$ satisfy
\[B(x)=e^{A(x)}.
\]
\end{proposition}

\noindent For the proof see \cite{Sta78,Sta98}.


\end{document}